\documentclass{m2an}
\usepackage{graphicx,fancyhdr}
\usepackage{graphics,color}
\usepackage{subfigure}
\usepackage{booktabs}
\usepackage{graphicx}
\graphicspath{{figure/}}  

\usepackage{color}
\usepackage{colortbl}
\usepackage{xcolor}
\definecolor{gray}{rgb}{0.85,0.85,0.85}
\usepackage[dvipdfm,linkcolor=blue,pdfstartview=FitH,
            CJKbookmarks=true,
            colorlinks=true, 
            pdfborder=000,   
            citecolor=blue]{hyperref}

\usepackage{syntonly}

\newtheorem{theorem}{Theorem}[section]

\newtheorem{lemma}{Lemma}[section]

\newtheorem{remark}{Remark}[section]

\numberwithin{equation}{section}
\begin{document}
\title{JACOBI-PREDICTOR-CORRECTOR APPROACH FOR THE FRACTIONAL ORDINARY DIFFERENTIAL EQUATIONS}
%
\author{Lijing Zhao}
\address{School of Mathematics and Statistics, Lanzhou University, Lanzhou 730000, China.\\
\textbf{zhaolj10@lzu.edu.cn}}
\author{Weihua Deng}
\address{School of Mathematics and Statistics, Lanzhou
University, Lanzhou 730000, China.\\
\textbf{dengwh@lzu.edu.cn}}
%
\date{}
\begin{abstract}
We present a novel numerical method, called {\tt
Jacobi-predictor-corrector approach}, for the numerical solution of
fractional ordinary differential equations based on the polynomial
interpolation and the Gauss-Lobatto quadrature w.r.t. the
Jacobi-weight function $\omega(s)=(1-s)^{\alpha-1}(1+s)^0$. This
method has the computational cost $O(N)$ and the convergent order
$IN$, where $N$ and $IN$ are, respectively, the total computational
steps and the number of used interpolating points. The detailed
error analysis is performed, and the extensive numerical experiments
confirm the theoretical results and show the robustness of this
method.
\end{abstract}
%
%
\maketitle
\section{Introduction}\label{sec:1}
This paper further discusses the numerical algorithm for the
following initial value problem
\begin{equation}\label{equa1.1}
D^\alpha_\ast x(t)=f\big(t,x(t)\big),~~~~~~x^{(k)}(0)=x^{(k)}_0,~~~k=0,1,\cdots,\lceil \alpha
\rceil-1,
\end{equation}
where $\alpha\in(0,\infty)$, $x_0^{(k)}$ can be any real numbers and $D_*^\alpha$ denotes the fractional derivative in the Caputo sense \cite{Podlubny:99},
defined by
\begin{equation*}
D_*^\alpha x(t) ={} _0D_t^{\alpha-n} D_t^n x(t)= \bigg\{
\begin{array}{ll}\frac{1}{\Gamma(n-\alpha)}\int_0 ^t
(t-\tau)^{n-\alpha-1}x^{(n)}(\tau)d\tau,& n-1 < \alpha< n;\\
\frac{d^{n}x(t)}{dt^{n}},      &\alpha=n;
\end{array}
\end{equation*}
where $n := \lceil \alpha \rceil$ is just the value $\alpha$ rounded up to the nearest integer, $D_t^n$ is the classical $n$th-order derivative, and
\begin{equation*}
{}_0D_t^{-\mu} x(t)=\frac{1}{\Gamma(\mu)}\int_0 ^t
(t-\tau)^{\mu-1}x(\tau)d\tau
\end{equation*}
is the Riemann-Liouville integral operator of order $\mu>0$. It is
well known that the initial value problem (\ref{equa1.1}) is
equivalent to the Volterra integral equation
\cite{Daftardar:04,Diethelm:02,Diethelm:04}
\begin{equation}\label{equa1.2}
x(t)=\sum_{k=0}^{\lceil \alpha \rceil-1}\frac{t^k}{k!}x_0^{(k)}+\frac{1}{\Gamma(\alpha)}\int^t_0(t-\tau)^{\alpha-1}f\big(\tau,x(\tau)\big)d\tau,~~~t\in[0,T],
\end{equation}
in the sense that if a continuous function solves (\ref{equa1.2}) if and only if it solves (\ref{equa1.1}).

Many approaches have been proposed to reslove (\ref{equa1.1}) or
(\ref{equa1.2}) numerically, such as
\cite{Ford:01,Diethelm:02,Diethelm:04,Deng:07,Deng:072}. Diethelm,
Ford and their coauthors successfully present the numerical
approximation of (\ref{equa1.2}) using Adams-type
predictor-corrector approach and give the corresponding detailed
error analysis in \cite{Diethelm:02} and \cite{Diethelm:04},
respectively. The convergent order of their approach is proved to be
$\min\{2,1+\alpha\}$, and the arithmetic complexity of their
algorithm with steps $N$ is $O(N^2)$, whereas a comparable algorithm
for a classical initial value problem only give rise to $O(N)$. The
challenge of the computational complexity is essentially because
fractional derivatives are non-local operators. This method has been
modified in \cite{Deng:07}, where the convergent order is improved
to be $\min\{2,1+2\alpha\}$ and almost half of the computational
cost is reduced, but the complexity is still $O(N^2)$. There are
already two typical ways which are suggested to overcome this
challenge. One seems to be the fixed memory principle of Podlubny
\cite{Podlubny:99}. However, it is shown that the fixed memory
principle is not suitable for Caputo derivative, because we cannot
reduce the computational cost significantly for preserving the
convergent order \cite{Diethelm:02,Ford:01}. The other more
promising idea seems to be the nested memory concept of Ford and
Simpson \cite{Ford:01,Deng:072} which can lead to
$O(N\log N)$ complexity, but still retain the order of
convergence. However, the convergent order there cannot exceed $2$.
For the effectiveness of the short memory principle, in
\cite{Ford:01}, $\alpha$ has to belong to the interval $(0,1)$; and
in \cite{Deng:072}, $\alpha$ must be within $(0,2)$.

In this work, we apprehend the Riemann-Liouville integral from the
viewpoint of a normal integral with a special weight function. Thus
we can deal with it based on the theories of the classical numerical
integration and  of polynomial interpolation \cite{Quarteroni:00}.
Then by using a predictor-corrector method, called {\tt
Jacobi-predictor-corrector approach}, we obtain a good numerical
approximation to (\ref{equa1.2}) with the convergent order $IN$,
which is the number of used interpolating points. Moreover, the
computational complexity is reduced to $O(N)$, the same as classical
initial value problem, which is one of the most exciting and
significant advantages of this algorithm.

The organization of this paper is as follows. In Section
\ref{sec:2}, we present the Jacobi-predictor-corrector approach and
its detailed algorithm. In Section \ref{sec:3}, the error analysis
of the numerical scheme is discussed in detail. The algorithm is
simply modified in Section \ref{sec:4} to deal with the extreme
cases. Two numerical examples are given in Section \ref{sec:5} to
confirm the theoretical results and to demonstrate the robustness
the algorithm. Concluding remarks are given in Section \ref{sec:6}.
Finally, some Jacobi-Gauss-Lobatto nodes and weights w.r.t. the
Jacobi-weight functions $(1-s)^{\alpha-1}(1+s)^0$ used in the
numerical experiments are listed in Appendix at the end of this
paper.

\section{Jacobi-predictor-corrector approach}\label{sec:2}

In this section we shall derive the fundamental algorithm for numerically solving the initial value
problems with Caputo derivative. It is the following transformation other than (\ref{equa1.2})
itself that underlies this algorithm: {\setlength\arraycolsep{2pt}
\begin{eqnarray}\label{equa2.1}
x(t)&=&\sum_{k=0}^{\lceil \alpha \rceil-1}\frac{t^k}{k!}x_0^{(k)}+\frac{1}{\Gamma(\alpha)}\int^{t}_0(t-\tau)^{\alpha-1}f\big(\tau,x(\tau)\big)d\tau{}
\nonumber\\
          &=&\sum_{k=0}^{\lceil \alpha \rceil-1}\frac{t^k}{k!}x_0^{(k)}+\frac{1}{\Gamma(\alpha)}\big(\frac{t}{2}\big)^\alpha\int^{1}_{-1}(1-s)^{\alpha-1}\tilde{f}\big(s,\tilde{x}(s)\big)ds,
\end{eqnarray}}
where
\begin{equation*}
\tilde{f}\big(s,\tilde{x}(s)\big)=f\Big(\frac{t}{2}(1+s),x\big(\frac{t}{2}(1+s)\big)\Big),~~-1\leq s \leq 1;
\end{equation*}
\begin{equation*}
\tilde{x}(s)=x\big(\frac{t}{2}(1+s)\big), ~~-1\leq s \leq 1.
\end{equation*}

We assume the function $f$ to be such that a unique solution of
(\ref{equa1.2}) exists in some interval $[0,T]$, specifically these
conditions are $(a)$ the continuity of $f$ with respect to both its
arguments and $(b)$ a Lipschitz condition with respect to the second
argument \cite{Diethelm:022}. Thus by the theory of the classical
numerical integration \cite{Quarteroni:00}, we can approximate the
integral in (\ref{equa2.1}) using the Jacobi-Gauss-Lobatto
quadrature w.r.t. the weight function
$\omega(s)=(1-s)^{\alpha-1}(1+s)^0$. That is,
\begin{equation}\label{equa2.2}
x(t)\approx \sum_{k=0}^{\lceil \alpha \rceil-1}\frac{t^k}{k!}x_0^{(k)}+\frac{1}{\Gamma(\alpha)}\big(\frac{t}{2}\big)^\alpha\sum_{j=0}^{JN}\omega_{j}\tilde{f}\big(s_j,\tilde{x}(s_j)\big),
\end{equation}
where $JN+1$, $\{\omega_{j}\}_{j=0}^{JN}$, $\{s_j\}_{j=0}^{JN}$ in
(\ref{equa2.2}) correspond to the number of, the weights of, and the
value of the Jacobi-Gauss-Lobatto nodes in the reference interval
$[-1,1]$, respectively.

Let us define a grid in the interval $[0,T]$ with $N+1$ equi-spaced
nodes $t_i$, given by
\begin{equation}\label{equa2.3}
t_i=ih,~~~i=0,\cdots,N,
\end{equation}
where $h=T/N$ is the step-length. Suppose that we have got the numerical values of $x(t)$ at
$t_0,t_1,\cdots,t_n$, which are denoted as $x_0,x_1,\cdots,x_n$, separately ($x_0=x_0^{(0)}$), now
we are going to compute the value of $x(t)$ at $t_{n+1}$, i.e. $x_{n+1}$.

By (\ref{equa2.2}),
\begin{equation}\label{27}
 x(t_{n+1})\approx \sum_{k=0}^{\lceil \alpha \rceil-1}\frac{t_{n+1}^k}{k!}x_0^{(k)}+\frac{1}{\Gamma(\alpha)}\big(\frac{t_{n+1}}{2}\big)^\alpha
                   \sum_{j=0}^{JN}\omega_{j}\tilde{f}_{n+1}\big(s_j,\tilde{x}_{n+1}(s_j)\big),
\end{equation}
where
\begin{equation*}
\tilde{f}_{n+1}\big(s,\tilde{x}_{n+1}(s)\big)=f\Big(\frac{t_{n+1}}{2}(1+s),x\big(\frac{t_{n+1}}{2}(1+s)\big)\Big),~~-1\leq s \leq 1;
\end{equation*}
\begin{equation*}
\tilde{x}_{n+1}(s)=x\big(\frac{t_{n+1}}{2}(1+s)\big), ~~-1\leq s \leq 1.
\end{equation*}

To do the summation of the second term of (\ref{27}), first we need
to evaluate the value of $f$ at the point $t_{n+1}$ since
$\tilde{f}_{n+1}\big(s_{JN},\tilde{x}_{n+1}(s_{JN})\big)=f\big(t_{n+1},x(t_{n+1})\big)$.
This value can be numerically got by using the interpolation of $f$
w.r.t. $t$ based on the known values of $f$ at the equi-spaced nodes
which are in the ``neighborhood" of $t_{n+1}$. For the other values
of $\tilde{f}_{n+1}\big(s_j,\tilde{x}_{n+1}(s_j)\big),\,0 \leq j
\leq JN-1$, we can also obtain them based on the interpolation of
$f$ at the equi-spaced nodes located in the ``neighborhood" of $s_j$
(should be $(1+s_j)t_{n+1}/2$ as to variable $t$). Denote $IN$ as
the number of equi-spaced nodes used for the interpolation. To start
this procedure, the values $x_0,x_1,\cdots,x_{IN-1}$ should be known
first, and should be accurate enough for not deteriorating the
accuracy of the algorithm.

Now we make it clear what the ``neighborhood" means. For predicting the value of $f$ at $t_{n+1}$,
we use the values of $f$ at $t_{n-IN+1},\,t_{n-IN+2},\,\cdots,\,t_{n-1},\,t_n$. For getting the
values of $\tilde{f}_{n+1}$ at $s_j,\,0 \leq j \leq JN-1$, the way to choose the ``neighborhood"
equi-spaced nodes is as follows: $(\textrm{i})$ try to make  $|ln-rn|$ as small as possible;
$(\textrm{ii})$ make $ln \ge rn$ if possible, where $ln$ and $rn$ are respectively the number of
equi-spaced points used on the left and right hand side of $s_j$ (should be $(1+s_j)t_{n+1}/2$ as
to variable $t$), obviously, $ln+rn=IN$.
So far, we arrive at the predictor-corrector formulas of
(\ref{equa2.1}) as
\begin{equation}\label{equa2.4}
x_{n+1}=\sum_{k=0}^{\lceil \alpha
\rceil-1}\frac{t_{n+1}^k}{k!}x_0^{(k)}+\frac{1}{\Gamma(\alpha)}\big(\frac{t_{n+1}}{2}\big)^\alpha\Big(\sum_{j=0}^{JN-1}\omega_{j}\tilde{f}_{n+1,j}
+\omega_{JN}f(t_{n+1},x^P_{n+1})\Big),
\end{equation}
and
\begin{equation}\label{equa2.5}
x^P_{n+1}=\sum_{k=0}^{\lceil \alpha
\rceil-1}\frac{t_{n+1}^k}{k!}x_0^{(k)}+\frac{1}{\Gamma(\alpha)}\big(\frac{t_{n+1}}{2}\big)^\alpha\sum_{j=0}^{JN}\omega_{j}\tilde{f}^P_{n+1,j},
\end{equation}
where $\{\tilde{f}^P_{n+1,j}\}_{j=0}^{JN}$ in (\ref{equa2.5}) means that all the values of
$\tilde{f}_{n+1}$ at the Jacobi-Gauss-Lobatto nodes are got by using the interpolations based on
the values of $\{f(t_i,x_i)\}_{i=0}^n$; whereas $\{\tilde{f}_{n+1,j}\}_{j=0}^{JN-1}$ in
(\ref{equa2.4}) are obtained by using the interpolations based on the values of
$\{f(t_i,x_i)\}_{i=0}^n$ and $f(t_{n+1},x^P_{n+1})$. The algorithm for realizing (\ref{equa2.4})
and (\ref{equa2.5}) is detailedly described as:
\\

$\textbf{Step~0.}$~\textbf{Some notations:}
{\setlength\arraycolsep{2pt}
\begin{eqnarray*}
\textrm{sum} &:=& \textrm{saving the value of the second summation
in the right hand of (\ref{equa2.5}) (or (\ref{equa2.4}))},
\nonumber\\
                  && \textrm{initially sum}=0;
\nonumber\\
Pl(t)        &:=& \textrm{the interpolating polynomial passing
through} ~(t_0,f(t_0,x_0)),(t_1,f(t_1,x_1)),\cdots,
\nonumber\\
             &&    ~(t_{IN-1},f(t_{IN-1},x_{IN-1}));
             \nonumber\\
Pr_{n+1}(t)  &:=& \textrm{the interpolating polynomial passing
through} ~(t_{n-IN+1},f(t_{n-IN+1},x_{n-IN+1})),
\nonumber\\
             &&  (t_{n-IN+2},f(t_{n-IN+2},x_{n-IN+2})),\cdots,(t_n,f(t_n,x_n));
             \nonumber\\
Cr_{n+1}(t)  &:=& \textrm{the interpolating polynomial passing
through} ~(t_{n-IN+2},f(t_{n-IN+2},x_{n-IN+2})),
\nonumber\\
             &&  (t_{n-IN+3},f(t_{n-IN+3},x_{n-IN+3})),\cdots,(t_{n+1},f(t_{n+1},x_{n+1}^P));
\nonumber\\
ln           &:=& \lceil IN/2 \rceil,~\textrm{to evaluate}
~\tilde{f}_{n+1}\big(s_j,\tilde{x}_{n+1}(s_j)\big),~\textrm{the
expected number of the interpolating}
\nonumber\\
             &&   \textrm{equi-spaced nodes on the left hand side of} ~s_j ;
\nonumber\\
rn           &:=& \lfloor IN/2 \rfloor,~\textrm{to evaluate}
~\tilde{f}_{n+1}\big(s_j,\tilde{x}_{n+1}(s_j)\big),~\textrm{the
expected number of the interpolating}
\nonumber\\
             &&   \textrm{equi-spaced nodes on the right hand side of} ~s_j;
\nonumber\\
le           &:=& \textrm{the number of the interpolating
equi-spaced nodes that can be used on the left of } ~s_j;
\nonumber\\
P_{n+1}(t)   &:=& \textrm{the interpolating polynomial passing
through} ~(t_{le-ln},f(t_{le-ln},x_{le-ln})),
\nonumber\\
             && (t_{le-ln+1}, f(t_{le-ln+1},x_{le-ln+1})),\cdots, (t_{le+rn-1},f(t_{le+rn-1},x_{le+rn-1})).
\end{eqnarray*}}

$\textbf{Step~1.}$~\textbf{To start the procedure:}

Compute $x_1,x_2,\cdots,x_{IN-1}$ by a single
step method (e.g., the Improved-Adams' methods in \cite{Deng:07}) with a sufficiently small
step-length $h_0$ such that $x_i,\,i=1,2,\cdots,IN-1,$ are accurate enough for not deteriorating
the accuracy of the method we are discussing.
\\

$\textbf{Step~2.}$~\textbf{To predict:}

\textbf{do}~$j=0,\cdots,JN$

~~~~~\textbf{if}~$le \leq ln$ (the number of the equi-spaced nodes
located on the left hand of $s_j$ (should be $(1+s_j)t_{n+1}/2$ as
to variable $t$) is equal to / less than what we expect)

~~~~~~~~$sum=sum+\omega_{j}Pl\big((1+s_j)t_{n+1}/2\big)$

~~~~~\textbf{else}~~~\textbf{if}~~~~$le+rn \geq n+1$ (the number of
the equi-spaced nodes located on the right hand of $s_j$ (should be
$(1+s_j)t_{n+1}/2$ as to variable $t$) is equal to / less than what
we expect)

~~~~~~~~~~~~~~~~~~~$sum=sum+\omega_{j}Pr_{n+1}\big((1+s_j)t_{n+1}/2\big)$

~~~~~~~~~~~~~\textbf{else}

~~~~~~~~~~~~~~~~~~~~$sum=sum+\omega_{j}P_{n+1}\big((1+s_j)t_{n+1}/2\big)$

\textbf{enddo}

$x^P_{n+1}=\sum_{k=0}^{\lceil \alpha
\rceil-1}\frac{t_{n+1}^k}{k!}x_0^{(k)}+\frac{1}{\Gamma(\alpha)}\big(\frac{t_{n+1}}{2}\big)^\alpha\cdot
sum$;
\\

$\textbf{Step~3.}$~\textbf{To correct:}

\textbf{do}~$j=0,\cdots,JN-1$

~~~~~\textbf{if}~$le \leq ln$

~~~~~~~~$sum=sum+\omega_{j}Pl\big((1+s_j)t_{n+1}/2\big)$

~~~~~\textbf{else}~~~\textbf{if}~~~~$le+rn \geq n+2$

~~~~~~~~~~~~~~~~~~~$sum=sum+\omega_{j}Cr_{n+1}\big((1+s_j)t_{n+1}/2\big)$

~~~~~~~~~~~~~\textbf{else}

~~~~~~~~~~~~~~~~~~~~$sum=sum+\omega_{j}P_{n+1}\big((1+s_j)t_{n+1}/2\big)$

\textbf{enddo}

$x_{n+1}=\sum_{k=0}^{\lceil \alpha
\rceil-1}\frac{t_{n+1}^k}{k!}x_0^{(k)}+\frac{1}{\Gamma(\alpha)}\big(\frac{t_{n+1}}{2}\big)^\alpha\cdot
\big(sum+\omega_{JN}f(t_{n+1},x^P_{n+1})\big)$.
\\

We call the above algorithm \emph{Jacobi-predictor-corrector
approach}.
\\

Although the description of this algorithm seems tedious, its operation is simple and mechanical.
It can be observed that for the computation of $x_{n+1}$, only changeless $2(JN+1)$ values are needed,
each of which should be interpolated by changeless $IN$ nearby values; whereas in \cite{Diethelm:02,Deng:07},
it should take $O(n+1)$ multiplications and divisions.
In other words, the computational cost here has no
relation with the variable $n+1$, just depends on the number of the
interpolating nodes $IN$ and the number of Jacobi-Gauss-Lobatto
nodes $JN+1$, so, to approximate $x(T)$, the total computational
cost is $O(N)$, comparing with $O(N^2)$ in
\cite{Diethelm:02,Deng:07} and $O(N\log N)$ in
\cite{Ford:01,Deng:072}, which is one of the most significant
advantages of this algorithm.
\section{Error Analysis}\label{sec:3}
First, we introduce four notations that will be used in the
following error analysis. The piecewise interpolating polynomial
based on the $IN$ nodes of $\big\{\big(t_i,f(t_i,x_i)\big)\big\}_{i=0}^{n}$
is denoted by $AI^{IN,P}\tilde{f}_{n+1}(s)$, where $-1 \le s \le 1$;
the one based on the $IN$ nodes of
$\big\{\big(t_i,f(t_i,x_i)\big)\big\}_{i=0}^{n}$ and
$\big(t_{n+1},f(t_{n+1},x_{n+1}^P)\big)$ is wrote as
$AI^{IN}\tilde{f}_{n+1}(s)$, where $-1 \le s \le 1$; the one based
on the $IN$ nodes of $\big\{\big(t_i,f(t_i,x(t_i))\big)\big\}_{i=0}^{n}$ is
signified by $PI^{IN,P}\tilde{f}_{n+1}(s)$, where $-1 \le s \le 1$;
the one based on the $IN$ nodes of
$\big\{\big(t_i,f(t_i,x(t_i))\big)\big\}_{i=0}^{n+1}$ is denoted as
$PI^{IN}\tilde{f}_{n+1}(s)$, where $-1 \le s \le 1$. Note that $x_i$
is the numerical solution and $x(t_i)$ is the exact solution.

Here authors state that the idea of the error analysis below is inspired by that in Kai Diethelm,  Neville J. Ford
and  Alan D. Freed's paper \cite{Diethelm:04}.

\subsection{Some preliminaries and a useful lemma}\label{sec:3.1}
Let $\omega^{\alpha,\beta}(x)=(1-x)^\alpha
(1+x)^\beta,\alpha>-1,\beta>-1$ be a Jacobi-weight function in the
usual sense. As illustrated in
\cite{Canuto:06,Guo:06,Guo:01,Guo:04,Shen:06,Wan:06,Hesthaven:07}, the set of
Jacobi polynomials $\{J_n^{\alpha,\beta}(x)\}_{n=0}^\infty$ forms a
complete $L_{\omega^{\alpha,\beta}}^2(-1,1)$-orthogonal system,
where  $L_{\omega^{\alpha,\beta}}^2(-1,1)$ is a weighted space
defined by
\begin{equation}\label{equa3.11}
L_{\omega^{\alpha,\beta}}^2(-1,1)=\{v:v ~\textrm{is measurable and} \parallel v \parallel_{\omega^{\alpha,\beta}}<\infty\},
\end{equation}
equipped with the norm
\begin{equation}
\parallel v \parallel_{\omega^{\alpha,\beta}}=\Big(\int_{-1}^{1}|v(x)|^2\omega^{\alpha,\beta}(x)dx\Big)^{\frac{1}{2}},
\end{equation}
and the inner product
\begin{equation}\label{equa3.12}
(u,v)_{\omega^{\alpha,\beta}}=\int_{-1}^{1}u(x)v(x)\omega^{\alpha,\beta}(x)dx.
\end{equation}

For bounding the approximation error of Jacobi polynomials, we need
the following non-uniformly-weighted Sobolev spaces as in
\cite{Shen:06}:
\begin{equation*}
H_{\omega^{\alpha,\beta},*}^m(-1,1):=\{v:\partial_x^k v \in L_{\omega^{\alpha+k,\beta+k}}^2(-1,1),~0\leq k\leq m\},
\end{equation*}
equipped with the inner product and the norm as
\begin{equation}\label{equa3.13}
(u,v)_{m,\omega^{\alpha,\beta},*}=\sum_{k=0}^m(\partial_x^k u,\partial_x^k v)_{\omega^{\alpha+k,\beta+k}},
\end{equation}
and
\begin{equation}\label{equa3.14}
\parallel v \parallel_{m,\omega^{\alpha,\beta},*}=\sqrt {(v,v)_{m,\omega^{\alpha,\beta},*}}.
\end{equation}
For any continuous functions $u$ and $v$ on $[-1,1]$, we define a
discrete inner product as
\begin{equation}\label{equa3.15}
(u,v)_N=\sum_{j=0}^Nu(x_j)v(x_j)\omega_j,
\end{equation}
where $\{\omega_j\}_{j=0}^N$ is a set of  Jacobi weights. The
following result follows from Lemma 3.3 in \cite{Chen:10}.
\begin{lemma}\label{lm3.11}
If $v\in H_{\omega^{\alpha,\beta},*}^m(-1,1)$ for some $m\geq 1$ and $\phi \in \mathcal{P}_N$,
then  for the Jacobi-Gauss-Lobatto integration, we have
\begin{equation*}
\mid(v,\phi)_{\omega^{\alpha,\beta}}-(v,\phi)_N\mid \leq CN^{-m}\parallel\partial_x^m v\parallel _{\omega^{\alpha+m-1,\beta+m-1}}\parallel\phi\parallel_{\omega^{\alpha,\beta}}.
\end{equation*}
\end{lemma}

\subsection{Auxiliary results}\label{sec:3.2}
By the definitions of the inner product (\ref{equa3.12}), the discrete inner product
(\ref{equa3.15}), and the notations given at the beginning of this section, we can rewrite
(\ref{equa2.1}) at $t=t_{n+1}$, (\ref{equa2.4}), and (\ref{equa2.5}), respectively, as
\begin{equation}\label{equa3.21}
x(t_{n+1})=\sum_{k=0}^{\lceil \alpha \rceil-1}\frac{t_{n+1}^k}{k!}x_0^{(k)}+\frac{1}{\Gamma(\alpha)}\big(\frac{t_{n+1}}{2}\big)^\alpha\bigg(\tilde{f}_{n+1}\big(\cdot,\tilde{x}_{n+1}(\cdot)\big),1\bigg)_{\omega^{\alpha-1,0}},
\end{equation}
\begin{equation}\label{equa3.22}
x_{n+1}=\sum_{k=0}^{\lceil \alpha \rceil-1}\frac{t_{n+1}^k}{k!}x_0^{(k)}+\frac{1}{\Gamma(\alpha)}\big(\frac{t_{n+1}}{2}\big)^\alpha\bigg(AI^{IN}\tilde{f}_{n+1}(\cdot),1\bigg)_{JN},
\end{equation}
and
\begin{equation}\label{equa3.23}
x^P_{n+1}=\sum_{k=0}^{\lceil \alpha \rceil-1}\frac{t_{n+1}^k}{k!}x_0^{(k)}+\frac{1}{\Gamma(\alpha)}\big(\frac{t_{n+1}}{2}\big)^\alpha\bigg(AI^{IN,P}\tilde{f}_{n+1}(\cdot),1\bigg)_{JN}.
\end{equation}

On the other hand, since each $\big\{ AI^{IN,P}\tilde{f}_{n+1}(s_j)
\big\}$ or $\big\{ AI^{IN}\tilde{f}_{n+1}(s_j) \big\}$ is
essentially a linear combination of parts of
$\big\{f(t_i,x_i)\big\}_{i=0}^{n}$ or
$\big\{f(t_i,x_i)\big\}_{i=0}^n$ and $f(t_{n+1},x_{n+1}^P)$, we can
also formally rewrite (\ref{equa2.4}) and (\ref{equa2.5}) as
\begin{equation}\label{equa3.24}
x_{n+1}=\sum_{k=0}^{\lceil \alpha
\rceil-1}\frac{t_{n+1}^k}{k!}x_0^{(k)}+\frac{1}{\Gamma(\alpha)}\big(\frac{t_{n+1}}{2}\big)^\alpha\Big[\sum_{i=0}^{n}a_{i,n+1}f(t_i,x_i)
+a_{n+1,n+1}f(t_{n+1},x^P_{n+1})\Big],
\end{equation}
and
\begin{equation}\label{equa3.25}
x^P_{n+1}=\sum_{k=0}^{\lceil \alpha
\rceil-1}\frac{t_{n+1}^k}{k!}x_0^{(k)}+\frac{1}{\Gamma(\alpha)}\big(\frac{t_{n+1}}{2}\big)^\alpha\sum_{i=0}^{n}b_{i,n+1}f(t_i,x_i),
\end{equation}
where $\{a_{i,n+1}\}_{i=0}^{n+1}$ and $\{b_{i,n+1}\}_{i=0}^{n}$ are
sets of real numbers depending on the number of the interpolating
nodes $IN$ and the positions of those Jacobi nodes in the interval
$[0,t_{n+1}]$. The formulae (\ref{equa3.24}) and (\ref{equa3.25})
can help us to understand the error analysis we will be performing.
First, we have the following estimate.
\\
\begin{theorem}\label{the3.21}
If ~$f(t,x)$ is sufficiently smooth on a suitable set $S\in
\mathbb{R}^2$, and $x(t)$ is also regular enough w.r.t. $t$,
then there is a constant $C_1$, independent of $n$ and $h$, such
that
\begin{equation}\label{equa3.29}
\left|
\frac{1}{\Gamma(\alpha)}\int_0^{t_{n+1}}(t_{n+1}-\tau)^{\alpha-1}f\big(\tau,x(\tau)\big)d\tau-\frac{1}
{\Gamma(\alpha)}\big(\frac{t_{n+1}}{2}\big)^\alpha\sum_{i=0}^{n}b_{i,n+1}f\big(t_i,x(t_i)\big)\right|
\leq C_1 t_{n+1}^\alpha h^{IN}.
\end{equation}
\end{theorem}
{\it Proof.} By the definitions of $\big\{ PI^{IN,P}\tilde{f}_{n+1}(s_j) \big\}$
and (\ref{equa3.15}), we have {\setlength\arraycolsep{2pt}
\begin{eqnarray}\label{equa3.26}
&&\frac{1}{\Gamma(\alpha)}\int_0^{t_{n+1}}(t_{n+1}-\tau)^{\alpha-1}f\big(\tau,x(\tau)\big)d\tau-\frac{1}{\Gamma(\alpha)}\big(\frac{t_{n+1}}{2}\big)^\alpha\sum_{i=0}^{n}b_{i,n+1}f(t_i,x_i)
\nonumber\\
&=&\frac{1}{\Gamma(\alpha)}\big(\frac{t_{n+1}}{2}\big)^\alpha
              \bigg[\bigg(\tilde{f}_{n+1}\big(\cdot,\tilde{x}_{n+1}(\cdot)\big),1\bigg)_{\omega^{\alpha-1,0}}-\bigg(PI^{IN,P}\tilde{f}_{n+1}(\cdot),1\bigg)_{JN}\bigg]{}
\nonumber\\
           &=&\frac{1}{\Gamma(\alpha)}\big(\frac{t_{n+1}}{2}\big)^\alpha
              \bigg[\bigg(\tilde{f}_{n+1}\big(\cdot,\tilde{x}_{n+1}(\cdot)\big),1\bigg)_{\omega^{\alpha-1,0}}-
              \bigg(\tilde{f}_{n+1}\big(\cdot,\tilde{x}_{n+1}(\cdot)\big),1\bigg)_{JN}
\nonumber\\
                  & & ~~~~~~~~~~~~~~~~~~~+\bigg(\tilde{f}_{n+1}\big(\cdot,\tilde{x}_{n+1}(\cdot)\big),1\bigg)_{JN}
                                          -\bigg(PI^{IN,P}\tilde{f}_{n+1}(\cdot),1\bigg)_{JN}\bigg]{}
\nonumber\\
           &=& I_{n+1,1}+I_{n+1,2},
\end{eqnarray}}
where
{\setlength\arraycolsep{2pt}
\begin{eqnarray}\label{equa3.27}
I_{n+1,1}&=&\frac{1}{\Gamma(\alpha)}\big(\frac{t_{n+1}}{2}\big)^\alpha
            \bigg[\bigg(\tilde{f}_{n+1}\big(\cdot,\tilde{x}_{n+1}(\cdot)\big),1\bigg)_{\omega^{\alpha-1,0}}-
              \bigg(\tilde{f}_{n+1}\big(\cdot,\tilde{x}_{n+1}(\cdot)\big),1\bigg)_{JN}\bigg],
\nonumber\\
I_{n+1,2}&=&\frac{1}{\Gamma(\alpha)}\big(\frac{t_{n+1}}{2}\big)^\alpha
            \bigg[\bigg(\tilde{f}_{n+1}\big(\cdot,\tilde{x}_{n+1}(\cdot)\big),1\bigg)_{JN}
                                          -\bigg(PI^{IN,P}\tilde{f}_{n+1}(\cdot),1\bigg)_{JN}\bigg].
\end{eqnarray}}

Under the assumption of that $f(t,x)$ is sufficiently smooth w.r.t.
$t$, from Lemma \ref{lm3.11}, we know $I_{n+1,1}$ can be
sufficiently small (because $m$ can be arbitrarily large), say, machine accuracy. Using the theories of the
Lagrange interpolation and of Gauss quadrature \cite{Quarteroni:00},
we have {\setlength\arraycolsep{2pt}
\begin{eqnarray}\label{equa3.28}
\mid I_{n+1,2}\mid &=& \frac{1}{\Gamma(\alpha)}\big(\frac{t_{n+1}}{2}\big)^\alpha
                      \left|\sum_{j=0}^{JN}\omega_j
                      \bigg(\tilde{f}_{n+1}\big(s_j,\tilde{x}_{n+1}(s_j)\big)-PI^{IN,P}\tilde{f}_{n+1}(s_j)\bigg)\right|
\nonumber\\
                    &\leq& \frac{1}{\Gamma(\alpha)}\big(\frac{t_{n+1}}{2}\big)^\alpha
                      \max_{0\leq j\leq JN}\mid\tilde{f}_{n+1}\big(s_j,\tilde{x}_{n+1}(s_j)\big)-PI^{IN,P}\tilde{f}_{n+1}(s_j)\mid\sum_{j=0}^{JN}\omega_j
\nonumber\\
\nonumber\\
                    &\leq&C(f,\alpha,IN,JN) t_{n+1}^\alpha h^{IN}.
\end{eqnarray}}
The inequalities in (\ref{equa3.28}) hold because of the
positiveness of the Gauss quadratures and
$\sum_{j=0}^{JN}\omega_j=\int_{-1}^1(1-s)^{\alpha-1}ds=2^{\alpha}/\alpha$.

Next we come to a result corresponding to the corrector formula.
Since the proof of the following theorem is very similar to the
above one, we omit it.
\begin{theorem}\label{the3.22}
If $f(t,x)$ is sufficiently smooth on a suitable set $S\in
\mathbb{R}^2$, and $x(t)$ is also regular enough w.r.t. $t$, then
there is a constant $C_2$, independent of $n$ and $h$, such that
\begin{equation}\label{equa3.210}
\Big|
\frac{1}{\Gamma(\alpha)}\int_0^{t_{n+1}}(t_{n+1}-\tau)^{\alpha-1}f\big(\tau,x(\tau)\big)d\tau-\frac{1}{\Gamma(\alpha)}
\big(\frac{t_{n+1}}{2}\big)^\alpha\sum_{i=0}^{n+1}a_{i,n+1}f\big(t_i,x(t_i)\big)\Big|
\leq C_2 t_{n+1}^\alpha h^{IN}.
\end{equation}
\end{theorem}

\subsection{Error analysis for the Jacobi-predictor-corrector approach}\label{sec:3.3}
In this subsection, we present the main result on the error estimate
of the Jacobi-predictor-corrector approach,
the proof of which is based on the mathematical induction and the
results given in the Subsection \ref{sec:3.2}. Through the following
result we can observe anther significant advantage of the presented
method---the convergence order is exactly the number of the
interpolating nodes $IN$, in other words, you can get your desired
convergent order just by choosing the enough number of interpolating
nodes.
\\

\begin{theorem}\label{the3.31}
If $f(t,x)$ is sufficiently smooth on a suitable set $S\in \mathbb{R}^2$, $h\leq1$, and $x(t)$ is regular
enough w.r.t. $t$, then for the Jacobi-predictor-corrector approach (\ref{equa2.4}) and
(\ref{equa2.5}) and for some suitably chosen $T$, there is a constant $C$ (depends on $\alpha,~IN$ and $JN$), independent of $n$ and
$h$, such that
\begin{equation}\label{equa3.31}
\max_{1\leq n+1\leq N}\mid x(t_{n+1})-x_{n+1}\mid \leq C h^{IN},
\end{equation}
where $N=T/h$.
\end{theorem}
{\it Proof.} We use the mathematical induction to prove this
theorem.

\textbf{a)} First we prove Eq. (\ref{equa3.31}) holds when $n+1=IN$:
Since $f$ is sufficiently smooth, $f$ is legitimately Lipschitz
continuous. Denoting $\tilde{L}$ as the Lipschitz constant of $f$
w.r.t. its second parameter $x$, then by (\ref{equa2.1}),
(\ref{equa3.25}) and Theorem \ref{the3.21}, there exists
{\setlength\arraycolsep{2pt}
\begin{eqnarray}\label{equa3.32}
&&\mid x(t_{n+1})-x^P_{n+1}\mid
\nonumber\\
&=&\frac{1}{\Gamma(\alpha)}\left|\int_0^{t_{n+1}}(t_{n+1}-\tau)^{\alpha-1}f\big(\tau,x(\tau)\big)d\tau-\big(\frac{t_{n+1}}{2}\big)^\alpha
\sum_{i=0}^{n}b_{i,n+1}f(t_i,x_i)\right|
\nonumber\\
&\leq&\frac{1}{\Gamma(\alpha)}\left|\int_0^{t_{n+1}}(t_{n+1}-\tau)^{\alpha-1}f\big(\tau,x(\tau)\big)d\tau-\big(\frac{t_{n+1}}{2}\big)^\alpha
\sum_{i=0}^{n}b_{i,n+1}f\big(t_i,x(t_i)\big)\right|
\nonumber\\
                   & &+\frac{1}{\Gamma(\alpha)}\big(\frac{t_{n+1}}{2}\big)^\alpha
                      \left|\sum_{j=0}^{JN}\omega_j
                      \bigg(PI^{IN,P}\tilde{f}_{n+1}(s_j)-AI^{IN,P}\tilde{f}_{n+1}(s_j)\bigg)\right|
\nonumber\\
                    &\leq& C_1 t_{n+1}^\alpha h^{IN}+\frac{1}{\Gamma(\alpha)}\big(\frac{t_{n+1}}{2}\big)^\alpha
                      \max_{0\leq j\leq JN}\mid PI^{IN,P}\tilde{f}_{n+1}(s_j)-AI^{IN,P}\tilde{f}_{n+1}(s_j)\mid\cdot\sum_{j=0}^{JN}\omega_j
\nonumber\\
                   &\leq& C_1 t_{n+1}^\alpha h^{IN}+\frac{t_{n+1}^\alpha}{\Gamma(\alpha+1)} C(IN,JN)\cdot\max_{0\leq i\leq n}\mid f\big(t_i,x(t_i)\big)-f(t_i,x_i)\mid\
\nonumber\\
                   &\leq& C_1 t_{n+1}^\alpha h^{IN}+\frac{\tilde{L}\cdot t_{n+1}^\alpha}{\Gamma(\alpha+1)} C(IN,JN)\cdot\max_{0\leq i\leq n}\mid x(t_i)-x_i\mid
\nonumber\\
                    &=&C_1 t_{n+1}^\alpha h^{IN}+C_3\tilde{L}t_{n+1}^\alpha\cdot\max_{0\leq i \leq n}\mid x(t_i)-x_i\mid.
\nonumber\\
\end{eqnarray}}

We have assumed that the starting error $\max_{0\leq i\leq
n=IN-1}\mid x(t_i)-x_i\mid $ is very small (not deteriorating the
accuracy of the present algorithm), so the first term in the right
hand of the last formula in (\ref{equa3.32}) plays the leading role,
thus
\begin{equation}\label{equa3.33}
\mid x(t_{n+1})-x^P_{n+1}\mid
\leq C_4 t_{n+1}^\alpha h^{IN}.
\end{equation}

Combining the above estimate with (\ref{equa2.1}), (\ref{equa3.24})
and Theorem \ref{the3.22}, {\setlength\arraycolsep{2pt}
\begin{eqnarray}\label{equa3.34}
&&\mid x(t_{n+1})-x_{n+1}\mid
\nonumber\\
&=&\frac{1}{\Gamma(\alpha)}\left|\int_0^{t_{n+1}}(t_{n+1}-\tau)^{\alpha-1}f\big(\tau,x(\tau)\big)d\tau-\big(\frac{t_{n+1}}{2}\big)^\alpha
\Big[\sum_{i=0}^{n}a_{i,n+1}f(t_i,x_i)+a_{n+1,n+1}f(t_{n+1},x^P_{n+1})\Big]\right|
\nonumber\\
&\leq&\frac{1}{\Gamma(\alpha)}\left|\int_0^{t_{n+1}}(t_{n+1}-\tau)^{\alpha-1}f\big(\tau,x(\tau)\big)d\tau-\big(\frac{t_{n+1}}{2}\big)^\alpha
\sum_{i=0}^{n+1}a_{i,n+1}f\big(t_i,x(t_i)\big)\right|
\nonumber\\
                   & &+\frac{1}{\Gamma(\alpha)}\big(\frac{t_{n+1}}{2}\big)^\alpha
                      \left|\sum_{j=0}^{JN}\omega_j
                      \bigg(PI^{IN}\tilde{f}_{n+1}(s_j)-AI^{IN}\tilde{f}_{n+1}(s_j)\bigg)\right|
\nonumber\\
                    &\leq& C_2 t_{n+1}^\alpha h^{IN}+\frac{1}{\Gamma(\alpha)}\big(\frac{t_{n+1}}{2}\big)^\alpha
                      \max_{0\leq j\leq JN}\mid PI^{IN}\tilde{f}_{n+1}(s_j)-AI^{IN}\tilde{f}_{n+1}(s_j)\mid\cdot\sum_{j=0}^{JN}\omega_j
\nonumber\\
                &\leq& C_2 t_{n+1}^\alpha h^{IN}+\frac{t_{n+1}^\alpha}{\Gamma(\alpha+1)} C(IN,JN)
\nonumber\\
                &&\cdot\max\{\max_{0\leq i\leq n}\mid f\big(t_i,x(t_i)\big)-f(t_i,x_i)\mid, \mid f\big(t_{n+1},x(t_{n+1})\big)-f(t_{n+1},x^P_{n+1})\mid\}
\nonumber\\
                &\leq& C_2 t_{n+1}^\alpha h^{IN}+\frac{\tilde{L}\cdot t_{n+1}^\alpha}{\Gamma(\alpha+1)} C(IN,JN)\cdot\max\{\max_{0\leq i\leq n}\mid x(t_i)-x_i\mid,\mid x(t_{n+1})-x^P_{n+1}\mid\}
\nonumber\\
                &\leq& C_2 t_{n+1}^\alpha h^{IN}+C_5\tilde{L} t_{n+1}^{\alpha}\cdot\max\{\max_{0\leq i\leq n}\mid x(t_i)-x_i\mid,C_4 t_{n+1}^\alpha h^{IN}\}
\nonumber\\
                &\leq& (C_2+C_6\tilde{L} t_{n+1}^{\alpha})t_{n+1}^{\alpha} h^{IN}
\nonumber\\
                &=& (C_2+C_6\tilde{L} t_{IN}^{\alpha})t_{IN}^{\alpha} h^{IN}
\nonumber\\
                &=& (C_2+C_6\tilde{L} \cdot IN^{\alpha}h^{\alpha})IN^{\alpha}h^{\alpha}\cdot h^{IN}
\nonumber\\
                &\leq& (C_2+C_6\tilde{L} \cdot IN^{\alpha})IN^{\alpha}\cdot h^{IN}:=Ch^{IN},
\end{eqnarray}}
where the last inequality holds since $h\leq1$.

\textbf{b)} We further prove Eq. (\ref{equa3.31}) holds for any
$n+1>IN$: Assume that $\max_{0\leq i\leq n+1}\mid
x(t_{n+1})-x_{n+1}\mid\leq C h^{IN}$, then we are going to prove
that $\mid x(t_{n+2})-x_{n+2}\mid\leq Ch^{IN}$. Since the
discussions are similar to the ones given in \textbf{a)}, we briefly
present them, {\setlength\arraycolsep{2pt}
\begin{eqnarray}\label{equa3.35}
&&\mid x(t_{n+2})-x^P_{n+2}\mid
\nonumber\\
    &\leq&C_1 t_{n+2}^\alpha h^{IN}+C_3\tilde{L}t_{n+2}^\alpha\cdot\max_{0\leq i \leq n+1}\mid x(t_i)-x_i\mid
\nonumber\\
    &\leq&C_1 T^\alpha h^{IN}+C_3\tilde{L}T^\alpha\cdot\max_{0\leq i \leq n+1}\mid x(t_i)-x_i\mid
\nonumber\\
    &\leq& (C_1+C_3C\tilde{L})T^\alpha h^{IN},
\end{eqnarray}}
and
{\setlength\arraycolsep{2pt}
\begin{eqnarray}\label{equa3.36}
&&\mid x(t_{n+2})-x_{n+2}\mid
\nonumber\\
                &\leq& C_2 t_{n+2}^\alpha h^{IN}+C_5\tilde{L} t_{n+2}^{\alpha}\cdot\max\{\max_{0\leq i\leq n+1}\mid x(t_i)-x_i\mid,\mid x(t_{n+2})-x^P_{n+2}\mid\}
\nonumber\\
                &\leq& C_2 T^\alpha h^{IN}+C_5\tilde{L} T^{\alpha}\cdot\max\{C h^{IN},(C_1+C_3C\tilde{L})T^\alpha
                h^{IN}\};
\end{eqnarray}}
by choosing $T$ sufficiently small, we can make sure that
$C_3C\tilde{L}T^{\alpha},C_1T^{\alpha},C_5C\tilde{L}T^{\alpha}$, as
well as $C_2T^\alpha$ are all bounded by $C/2$, thus
{\setlength\arraycolsep{2pt}
\begin{eqnarray}\label{equa3.37}
&&\mid x(t_{n+2})-x_{n+2}\mid
\nonumber\\
                &\leq& C_2 T^\alpha h^{IN}+C_5\tilde{L}T^{\alpha}\cdot C  h^{IN}
\nonumber\\
                &\leq&Ch^{IN}.
\end{eqnarray}}

\begin{remark}\label{remark3.1}
In practical computations, the Improved-Adams' methods proposed in
\cite{Deng:07} can be used to start the algorithm. We can take the
step-length $h_0$ discussed in the algorithm description in Section
\ref{sec:2} as $h\cdot 10^{-k}$, where $k \geq 1$ is a given
integer, then by the result in \cite{Deng:07}, there exists
{\setlength\arraycolsep{2pt}
\begin{eqnarray}\label{equa3.38}
&&\max_{0\leq i \leq IN-1}\mid x(t_i)-x_i \mid
 \nonumber\\
&=&  O\left(h_0^{\min\{1+2\alpha,2\}}\right)
\nonumber\\
&=& \left\{ \begin{array}{ll}
10^{-k(1+2\alpha)}\cdot O(h^{1+2\alpha}),  &   \textrm{if}~~0< \alpha\leq0.5; \\
10^{-2k}\cdot O(h^2), &   \textrm{if}~~\alpha> 0.5.
\end{array} \right.
\end{eqnarray}}

If taking $h=10^{-m}$, where $m$ is a given positive integer, by a
simple computation, we obtain that $\mid
x(t_{n+1})-x_{n+1}\mid=O(h^{IN})$,
 as long as the integers $k$ and $m$ satisfy
\begin{equation}\label{equa3.9}
IN <\left\{ \begin{array}{ll}
(1+2\alpha)(1+\frac{k}{m}),  &   \textrm{if}~0< \alpha\leq0.5; \\
2+\frac{2k}{m}, &   \textrm{if}~~\alpha> 0.5.
\end{array} \right.
\end{equation}

\end{remark}

\section{Modifications of the Algorithm}\label{sec:4}

We have completed the description of the basic algorithm and its
most important properties. The convergent order of the algorithm is
exactly equal to the number of the interpolating points. As is well
known, in practice, the effectiveness of the algorithm is closely
related to the regularity of the solution of the equation. For the
fractional ordinary differential equation, most of the time, the
smoothness of its solution at $t=0$ is weaker than other places
\cite{Deng:10}. So we will simply discuss this issue in the
following. Another issue we will also mention is how to use this
algorithm when $\alpha$ is very close to $0$.

\subsection{The function $f$ or $x$ is not very smooth at the starting point}\label{sec:4.1}
When the smoothness of $f$ or $x$ is weaker at the initial time zero
than other time, the sensible way is to divide the interval $[0,T]$
into two parts $[0,T_0]$ and $[T_0,T]$, where $T_0$ is a small
positive real number. For the small interval $[0,T_0]$, we use the
Gauss-Lobatto quadrature with the weight function $\omega(s)=1$. For
the remaining part $[T_0,T]$, the algorithm provided above is
employed. That is,
{\setlength\arraycolsep{2pt}
\begin{eqnarray}\label{equa4.11}
x(t_{n+1})&=&\sum_{k=0}^{\lceil \alpha \rceil-1}\frac{t_{n+1}^k}{k!}x_0^{(k)}+\frac{1}{\Gamma(\alpha)}\int^{T_0}_0(t_{n+1}-\tau)^{\alpha-1}f\big(\tau,x(\tau)\big)d\tau
\nonumber\\
           &&+\frac{1}{\Gamma(\alpha)}\int^{t_{n+1}}_{T_0}(t_{n+1}-\tau)^{\alpha-1}f\big(\tau,x(\tau)\big)d\tau
\nonumber\\
          &=&\sum_{k=0}^{\lceil \alpha \rceil-1}\frac{t_{n+1}^k}{k!}x_0^{(k)}
            +\frac{1}{\Gamma(\alpha)}\int^{T_0}_0(t_{n+1}-\tau)^{\alpha-1}f\big(\tau,x(\tau)\big)d\tau
\nonumber\\
           &&+\frac{1}{\Gamma(\alpha)}\big(\frac{t_{n+1}-T_0}{2}\big)^\alpha\int^{1}_{-1}(1-s)^{\alpha-1}\tilde{f}_{n+1}\big(s,\tilde{x}_{n+1}(s)\big)ds
\nonumber\\
          &\approx&\sum_{k=0}^{\lceil \alpha \rceil-1}\frac{t_{n+1}^k}{k!}x_0^{(k)}
          +\frac{1}{\Gamma(\alpha)}\sum_{j=0}^{\tilde{JN}}\tilde{\omega}_j(t_{n+1}-\tau_j)^{\alpha-1}f\big(\tau_j,x(\tau_j)\big)
\nonumber\\
           &&+\frac{1}{\Gamma(\alpha)}\big(\frac{t_{n+1}-T_0}{2}\big)^\alpha\sum_{j=0}^{JN}\omega_j\tilde{f}_{n+1}\big(s_j,\tilde{x}_{n+1}(s_j)\big),
\end{eqnarray}}
where $\tilde{JN},~\{\tilde{\omega}_{j}\}_{j=0}^{\tilde{JN}}$ and $\{\tau_j\}_{j=0}^{\tilde{JN}}$
correspond to the number of, the weights of, and the values of the Gauss-Lobatto nodes with the
weight $\omega(s)=1$ in the interval $[0,T_0]$, respectively. The values of
$\big\{f\big(\tau_j,x(\tau_j)\big)\big\}_{j=0}^{\tilde{JN}}$ can be computed as in the starting
procedure. Since $f$ and $x$ are continuous in the interval $[0,T_0]$, by the theory of Gauss
quadrature \cite{Quarteroni:00} and the analysis above, we can see that if $\tilde{JN}$ is a big
number the accuracy of the total error still can be remained.

\subsection{The value of $\alpha$ is very small}\label{sec:4.2}
In this subsection, we discuss the case that $\alpha$ is very small, although it happens very less
often. When $\alpha$ becomes bigger, the weight of the Gauss-Lobatto quadrature at the endpoint of
the righthand side becomes smaller, the provided algorithm becomes more robust. When $\alpha$ is
very small, say, smaller than $0.1$, the weight at the endpoint of the righthand side of the
interval is much bigger than other places (see the Appendix),
which may impacts the
robustness of the algorithm. There are two choices to deal with this
problem: one is to try to avoid using the high order interpolation
in the algorithm; one is to divide the interval $[0,T]$ into two
subintervals $[0,T_1]$ and $[T_1,T]$, then similarly do what we
do in the last subsection, that is,
{\setlength\arraycolsep{2pt}
\begin{eqnarray}\label{equa4.21}
x(t)&=&\sum_{k=0}^{\lceil \alpha \rceil-1}\frac{t^k}{k!}x_0^{(k)}+\frac{1}{\Gamma(\alpha)}\int^{T_1}_0(t-\tau)^{\alpha-1}f\big(\tau,x(\tau)\big)d\tau
\nonumber\\
           &&+\frac{1}{\Gamma(\alpha)}\int^{t}_{T_1}(t-\tau)^{\alpha-1}f\big(\tau,x(\tau)\big)d\tau
\nonumber\\
          &=&\sum_{k=0}^{\lceil \alpha \rceil-1}\frac{t^k}{k!}x_0^{(k)}
            +\frac{1}{\Gamma(\alpha)}\int^{T_1}_0(t-\tau)^{\alpha-1}f\big(\tau,x(\tau)\big)d\tau
\nonumber\\
           &&+\frac{1}{\Gamma(\alpha)}\left(\frac{t-T_1}{2}\right)^\alpha\int^{1}_{-1}(1-s)^{\alpha-1}\tilde{f}\big(s,\tilde{x}(s)\big)ds
\nonumber\\
          &\approx&\sum_{k=0}^{\lceil \alpha \rceil-1}\frac{t^k}{k!}x_0^{(k)}
          +\frac{1}{\Gamma(\alpha)}\sum_{j=0}^{\tilde{JN}}\tilde{\omega}_j(t-\tau_j)^{\alpha-1}f\big(\tau_j,x(\tau_j)\big)
\nonumber\\
           &&+\frac{1}{\Gamma(\alpha)}\left(\frac{t-T_1}{2}\right)^\alpha\sum_{j=0}^{JN}\omega_j\tilde{f}\big(s_j,\tilde{x}(s_j)\big),
\end{eqnarray}}
where $\tilde{JN},~\{\tilde{\omega}_{j}\}_{j=0}^{\tilde{JN}}$ and $\{\tau_j\}_{j=0}^{\tilde{JN}}$
are the same as those defined in the last subsection. And
$\big\{f\big(\tau_j,x(\tau_j)\big)\big\}_{j=0}^{\tilde{JN}}$ in the second term of the righthand
side of the last formula can be computed by interpolation.

\section{Numerical Experiments}\label{sec:5}

In this section we present two numerical examples to verify the
convergent orders derived above and to demonstrate the robustness of
the provided methods for big $T$. We only consider the examples with
$0 <\alpha<2$, since the algorithm will be more robust for
$\alpha\geq2$. All numerical computations were done in Matlab 7.5.0 on a
normal laptop with 1GB of memory.

\subsection{Verification of the error estimates}\label{sec:5.1}
First, we use the following example to verify the convergent order:
\begin{equation}\label{equa5.11}
D_*^\alpha
x(t)=-x(t)+\frac{\Gamma(9)}{\Gamma{(9-\alpha)}}t^{8-\alpha}+3\cdot\frac{\Gamma(8)}{\Gamma{(8-\alpha)}}t^{7-\alpha}+t^{8}-3t^7,
\end{equation}
with the initial condition(s) $x(0)=0$ (and $x'(0)=0$ if $1\leq
\alpha <2$). The exact solution of this initial value problem is
\begin{equation}\label{equa5.12}
x(t)=t^{8}+3t^7.
\end{equation}

We start the procedure with the Improved-Adams' methods in \cite{Deng:07} as discussed in Remark
\ref{remark3.1}, i.e. the values of $x(t)$ at $t_0,t_1,\cdots,t_{IN-1}$ are computed by the
Improved-Adams' methods. The convergent orders are verified at $T=1$, and the number of the Jacobi
nodes is taken as $JN+1=27$. The number of the interpolating nodes $IN$ is respectively taken as
$2,\,3,\,4$ and $5$ to expect that the corresponding convergent order is also $2,\,3,\,4$ and $5$.
The numerical results of the maximum errors for the Jacobi-predictor-corrector approach are showed
in the following tables, where ``CO" means the convergent order. The nodes and weights used in the
Gauss-Lobatto quadrature w.r.t. the weight functions $\omega(s)=(1-s)^{\alpha-1}(1+s)^0$ for
various $\alpha$ are given in Appendix. Form the results in Table \ref{table5.11} to Table
\ref{table5.14}, we can see that the data confirm the theoretical results provided in Section
\ref{sec:3.3}.

Table \ref{table5.15} and Table \ref{table5.16} give the numerical results of the maximum errors
for the fractional Adams methods in \cite{Diethelm:02} and for the Improved Adams methods in \cite{Deng:07}.
The compare of Table \ref{table5.11} to Table \ref{table5.14}, Table \ref{table5.15} and Table \ref{table5.16}
indicates that the Jacobi-predictor-corrector approach is effective.

See the results in Table \ref{table5.13} and Table \ref{table5.14}
for $\alpha=0.1$, $IN=4$ and $5$; like what we discussed in Section
\ref{sec:4.2}, the results tell us that we must be more careful to
use the provided algorithm when $\alpha$ is very small (letting $IN$
be small or dividing the original interval into subintervals).
However, we still confirm the convergent order by taking small $T$
($T=0.1$ and $T=0.001$) in Table \ref{table5.17}.

\begin{table}[!htb]\fontsize{9.5pt}{12pt}\selectfont
\centering \caption{The maximum errors for (\ref{equa5.11}) when
$t\in[0,1]$ and $IN=2$.}\vspace{5pt}

\begin{tabular*}{\linewidth}{@{\extracolsep{\fill}}*{9}{c|cc|cc|cc|cc}}
 \toprule  
h & $\alpha=0.1$ & CO & $\alpha=0.3$ & CO & $\alpha=0.5$ & CO & $\alpha=0.7$ & CO \\
\midrule
1/10&    1.34 1e+0 & -     & 4.21 1e-1 & -    & 2.27 1e-1  & -    & 1.68 1e-1 & - \\
1/20&    3.75 1e-1 & 1.84  & 8.57 1e-2 & 2.30 & 4.65 1e-2  & 2.29 & 3.99 1e-2 & 2.07 \\
1/40&    8.40 1e-2 & 2.16  & 1.70 1e-2 & 2.34 & 9.99 1e-3  & 2.22 & 9.14 1e-3 & 2.13 \\
1/80&    1.71 1e-2 & 2.30  & 2.98 1e-3 & 2.51 & 1.89 1e-3  & 2.40 & 2.14 1e-3 & 2.10 \\
1/160&   3.74 1e-3 & 2.19  & 6.67 1e-4 & 2.16 & 4.17 1e-4  & 2.18 & 5.84 1e-4 & 1.87 \\
1/320&   6.97 1e-4 & 2.42  & 9.84 1e-5 & 2.76 & 1.06 1e-4  & 1.98 & 1.28 1e-4 & 2.19 \\
1/640&   1.49 1e-4 & 2.22  & 2.75 1e-5 & 1.84 & 2.78 1e-5  & 1.93 & 3.19 1e-5 & 2.00 \\
1/1280&  3.67 1e-5 & 2.03  & 7.84 1e-6 & 1.81 & 8.13 1e-6  & 1.77 & 9.03 1e-6 & 1.82 \\
1/2560&  9.30 1e-6 & 1.98  & 2.08 1e-6 & 1.91 & 1.94 1e-6  & 2.07 & 2.40 1e-6 & 1.91 \\
\midrule
\midrule
h & $\alpha=0.9$ & CO & $\alpha=1.2$ & CO & $\alpha=1.5$ & CO & $\alpha=1.8$ & CO \\
\midrule
1/10&    1.51 1e-1 & -     & 1.47 1e-1 & -    & 1.47 1e-1  & -    & 1.49 1e-1 & - \\
1/20&    4.03 1e-2 & 1.91  & 4.16 1e-2 & 1.83 & 3.95 1e-2  & 1.90 & 3.75 1e-2 & 1.99 \\
1/40&    8.84 1e-3 & 2.19  & 8.26 1e-3 & 2.33 & 9.01 1e-3  & 2.13 & 1.08 1e-2 & 1.80 \\
1/80&    2.32 1e-3 & 1.93  & 2.05 1e-3 & 2.01 & 2.35 1e-3  & 1.94 & 2.58 1e-3 & 2.06 \\
1/160&   6.16 1e-4 & 1.91  & 5.07 1e-4 & 2.02 & 6.14 1e-4  & 1.94 & 6.67 1e-4 & 1.95 \\
1/320&   1.48 1e-4 & 2.05  & 1.44 1e-4 & 1.82 & 1.64 1e-4  & 1.91 & 1.67 1e-4 & 1.99 \\
1/640&   3.84 1e-5 & 1.95  & 3.91 1e-5 & 1.88 & 3.77 1e-5  & 2.12 & 4.21 1e-5 & 1.99 \\
1/1280&  9.75 1e-6 & 1.98  & 1.09 1e-5 & 1.84 & 1.01 1e-5  & 1.90 & 1.11 1e-5 & 1.93 \\
1/2560&  2.46 1e-6 & 1.99  & 2.69 1e-6 & 2.02 & 2.73 1e-6  & 1.89 & 2.86 1e-6 & 1.95 \\
\bottomrule 
\end{tabular*}\label{table5.11}
\end{table}

\begin{table}[ht]\fontsize{9.5pt}{12pt}\selectfont
\centering \caption{The maximum errors for (\ref{equa5.11}) when
$t\in[0,1]$ and $IN=3$.}\vspace{5pt}
\begin{tabular*}{\linewidth}{@{\extracolsep{\fill}}*{9}{c|cc|cc|cc|cc}}
 \toprule  
h & $\alpha=0.1$ & CO & $\alpha=0.3$ & CO & $\alpha=0.5$ & CO & $\alpha=0.7$ & CO \\
\midrule
1/10&    6.47 1e-1 & -     & 1.50 1e-1 & -    & 6.69 1e-2  & -    & 4.26 1e-2 & - \\
1/20&    8.54 1e-2 & 2.92  & 1.59 1e-2 & 3.24 & 7.01 1e-3  & 3.25 & 5.23 1e-3 & 3.03 \\
1/40&    9.64 1e-3 & 3.15  & 1.57 1e-3 & 3.34 & 7.19 1e-4  & 3.28 & 5.72 1e-4 & 3.19 \\
1/80&    1.00 1e-3 & 3.27  & 1.37 1e-4 & 3.52 & 6.85 1e-5  & 3.39 & 7.19 1e-5 & 2.99 \\
1/160&   1.07 1e-4 & 3.22  & 1.39 1e-5 & 3.30 & 7.05 1e-6  & 3.28 & 9.77 1e-6 & 2.88 \\
1/320&   1.02 1e-5 & 3.39  & 1.06 1e-6 & 3.71 & 9.50 1e-7  & 2.89 & 1.05 1e-6 & 3.22 \\
1/640&   1.17 1e-6 & 3.13  & 1.45 1e-7 & 2.87 & 1.25 1e-7  & 2.93 & 1.36 1e-7 & 2.95 \\
1/1280&  1.45 1e-7 & 3.00  & 1.95 1e-8 & 2.90 & 1.76 1e-8  & 2.82 & 1.98 1e-8 & 2.78 \\
1/2560&  1.83 1e-8 & 2.99  & 2.57 1e-9 & 2.92 & 2.17 1e-9  & 3.02 & 2.45 1e-9 & 3.01 \\
\midrule
\midrule
h & $\alpha=0.9$ & CO & $\alpha=1.2$ & CO & $\alpha=1.5$ & CO & $\alpha=1.8$ & CO \\
\midrule
1/10&    3.51 1e-2 & -     & 3.27 1e-2 & -    & 3.24 1e-2  & -    & 3.27 1e-2 & - \\
1/20&    4.94 1e-3 & 2.83  & 4.84 1e-3 & 2.76 & 4.46 1e-3  & 2.86 & 4.24 1e-3 & 2.95 \\
1/40&    5.40 1e-4 & 3.19  & 5.29 1e-4 & 3.19 & 5.89 1e-4  & 2.92 & 6.79 1e-4 & 2.64 \\
1/80&    7.32 1e-5 & 2.88  & 6.70 1e-5 & 2.98 & 7.95 1e-5  & 2.89 & 8.04 1e-5 & 3.08 \\
1/160&   9.71 1e-6 & 2.91  & 9.10 1e-6 & 2.88 & 1.05 1e-5  & 2.92 & 1.09 1e-5 & 2.88 \\
1/320&   1.23 1e-6 & 2.99  & 1.19 1e-6 & 2.93 & 1.37 1e-6  & 2.94 & 1.38 1e-6 & 2.98 \\
1/640&   1.51 1e-7 & 3.02  & 1.65 1e-7 & 2.86 & 1.52 1e-7  & 3.17 & 1.68 1e-7 & 3.03 \\
1/1280&  2.08 1e-8 & 2.86  & 2.14 1e-8 & 2.95 & 2.11 1e-8  & 2.85 & 2.18 1e-8 & 2.95 \\
1/2560&  2.49 1e-9 & 3.07  & 2.69 1e-9 & 2.99 & 2.78 1e-9  & 2.92 & 3.04 1e-9 & 2.84 \\
\bottomrule 
\end{tabular*}\label{table5.12}
\end{table}
\begin{table}[ht]\fontsize{9.5pt}{12pt}\selectfont
\centering \caption{The maximum errors for (\ref{equa5.11}) when
$t\in[0,1]$ and $IN=4$.}\vspace{5pt}
\begin{tabular*}{\linewidth}{@{\extracolsep{\fill}}*{9}{c|cc|cc|cc|cc}}
 \toprule  
h & $\alpha=0.1$ & CO & $\alpha=0.3$ & CO & $\alpha=0.5$ & CO & $\alpha=0.7$ & CO \\
\midrule
1/10&    2.39 1e-1 & -     & 4.71 1e-2 & -    & 1.48 1e-2  & -    & 4.50 1e-3 & - \\
1/20&    1.61 1e-2 & 3.90  & 2.41 1e-3 & 4.28 & 5.09 1e-4  & 4.86 & 8.28 1e-5 & 5.76 \\
1/40&    1.08 1e-3 & 3.89  & 1.05 1e-4 & 4.52 & 1.43 1e-5  & 5.16 & 8.79 1e-6 & 3.24 \\
1/80&    3.05 1e-4 & 1.83  & 4.22 1e-6 & 4.64 & 2.99 1e-7  & 5.57 & 8.77 1e-7 & 3.33 \\
1/160&   4.30 1e-3 &-3.82  & 1.80 1e-7 & 4.55 & 1.73 1e-8  & 4.11 & 7.20 1e-8 & 3.61 \\
1/320&   6.31 1e-1 &-7.20  & 4.97 1e-9 & 5.18 & 3.92 1e-9  & 2.15 & 5.10 1e-9 & 3.82 \\
1/640&   1.94 1e+1 &-4.95  & 2.87 1e-10& 4.11 & 2.38 1e-10 & 4.04 & 3.36 1e-10& 3.92 \\
1/1280&  2.27 1e+4 &-10.2  & 1.57 1e-11& 4.19 & 1.93 1e-11 & 3.62 & 2.26 1e-11& 3.89 \\
1/2560&  1.15 1e+12&-25.6  & 1.04 1e-12& 3.92 & 1.03 1e-12 & 4.23 & 1.47 1e-12& 3.94 \\
\midrule
\midrule
h & $\alpha=0.9$ & CO & $\alpha=1.2$ & CO & $\alpha=1.5$ & CO & $\alpha=1.8$ & CO \\
\midrule
1/10&    7.94 1e-4 & -     & 2.82 1e-3 & -    & 4.34 1e-3  & -    & 5.00 1e-3 & - \\
1/20&    2.02 1e-4 & 1.97  & 3.29 1e-4 & 3.10 & 3.65 1e-4  & 3.57 & 3.80 1e-4 & 3.72 \\
1/40&    1.59 1e-5 & 3.67  & 2.06 1e-5 & 4.00 & 2.20 1e-5  & 4.06 & 2.69 1e-5 & 3.82 \\
1/80&    1.22 1e-6 & 3.70  & 1.31 1e-6 & 3.97 & 1.54 1e-6  & 3.83 & 1.57 1e-6 & 4.10 \\
1/160&   9.00 1e-8 & 3.76  & 8.89 1e-8 & 3.89 & 1.08 1e-7  & 3.84 & 1.05 1e-7 & 3.90 \\
1/320&   5.74 1e-9 & 3.97  & 5.58 1e-9 & 3.99 & 6.62 1e-9  & 4.03 & 6.78 1e-9 & 3.95 \\
1/640&   3.85 1e-10& 3.90  & 3.94 1e-10& 3.82 & 3.86 1e-10 & 4.10 & 4.35 1e-10& 3.96 \\
1/1280&  2.36 1e-11& 4.03  & 2.78 1e-11& 3.82 & 2.58 1e-11 & 3.90 & 2.70 1e-11& 4.01 \\
1/2560&  1.50 1e-12& 3.98  & 1.70 1e-12& 4.04 & 1.67 1e-12 & 3.95 & 1.80 1e-12& 3.91 \\
\bottomrule 
\end{tabular*}\label{table5.13}
\end{table}
\begin{table}[ht]\fontsize{9.5pt}{12pt}\selectfont
\centering \caption{The maximum errors for (\ref{equa5.11}) when
$t\in[0,1]$ and $IN=5$.}\vspace{5pt}
\begin{tabular*}{\linewidth}{@{\extracolsep{\fill}}*{9}{c|cc|cc|cc|cc}}
\toprule  
h & $\alpha=0.1$ & CO & $\alpha=0.3$ & CO & $\alpha=0.5$ & CO & $\alpha=0.7$ & CO \\
\midrule
1/10&    4.79 1e-2 & -     & 1.49 1e-2 & -    & 4.07 1e-3  & -    & 1.23 1e-3 & - \\
1/20&    1.51 1e-2 & 1.67  & 3.54 1e-4 & 5.40 & 7.22 1e-5  & 5.82 & 1.49 1e-5 & 6.37 \\
1/40&    8.11 1e-1 &-5.75  & 7.74 1e-6 & 5.52 & 1.08 1e-6  & 6.07 & 2.74 1e-7 & 5.76 \\
1/80&    1.25 1e+4 &-13.9  & 1.58 1e-7 & 5.61 & 1.39 1e-8  & 6.28 & 1.68 1e-8 & 4.03 \\
1/160&   -         &-      & 3.31 1e-9 & 5.58 & 1.93 1e-10 & 6.17 & 7.47 1e-10& 4.49 \\
1/320&   -         &-      & 5.19 1e-11& 5.99 & 1.95 1e-11 & 3.31 & 2.57 1e-11& 4.86 \\
1/640&   -         &-      & 1.60 1e-12& 5.02 & 5.64 1e-13 & 5.11 & 8.92 1e-13& 4.85 \\
1/1280&  -         &-      & 5.06 1e-14& 4.98 & 1.87 1e-14 & 4.92 & 3.06 1e-14& 4.86 \\
\midrule
\midrule
h & $\alpha=0.9$ & CO & $\alpha=1.2$ & CO & $\alpha=1.5$ & CO & $\alpha=1.8$ & CO \\
\midrule
1/10&    2.69 1e-4 & -     & 4.47 1e-4 & -    & 7.36 1e-4  & -    & 8.87 1e-4  & - \\
1/20&    1.49 1e-5 & 4.18  & 2.60 1e-5 & 4.10 & 2.94 1e-5  & 4.64 & 3.14 1e-5  & 4.82 \\
1/40&    6.52 1e-7 & 4.51  & 8.53 1e-7 & 4.93 & 9.18 1e-7  & 5.00 & 1.10 1e-6  & 4.84 \\
1/80&    2.47 1e-8 & 4.72  & 2.63 1e-8 & 5.02 & 3.25 1e-8  & 4.82 & 3.07 1e-8  & 5.16 \\
1/160&   9.01 1e-10& 4.78  & 9.88 1e-10& 4.73 & 1.14 1e-9  & 4.84 & 1.08 1e-9  & 4.83 \\
1/320&   2.90 1e-11& 4.96  & 2.86 1e-11& 5.11 & 3.31 1e-11 & 5.10 & 3.42 1e-11 & 4.98 \\
1/640&   9.53 1e-13& 4.93  & 1.01 1e-12& 4.83 & 1.00 1e-12 & 5.05 & 1.10 1e-12 & 4.96\\
1/1280&  3.46 1e-14& 4.78  & 3.38 1e-14& 4.90 & 3.24 1e-14 & 4.95 & 3.46 1e-14 & 4.99 \\
\bottomrule 
\end{tabular*}\label{table5.14}
\end{table}

\begin{table}[ht]\fontsize{9.5pt}{12pt}\selectfont
\centering \caption{The maximum errors of fractional Adams methods for (\ref{equa5.11}) when
$t\in[0,1]$.}\vspace{5pt}
\begin{tabular*}{\linewidth}{@{\extracolsep{\fill}}*{9}{c|cc|cc|cc|cc}}
\toprule  
h & $\alpha=0.1$ & CO & $\alpha=0.3$ & CO & $\alpha=0.5$ & CO & $\alpha=0.7$ & CO \\
\midrule
1/10&    2.09 1e+0 & -     & 8.45 1e-1 & -    & 4.51 1e-1  & -    & 2.89 1e-1 & - \\
1/20&    1.17 1e+0 & 0.83  & 3.32 1e-1 & 1.35 & 1.46 1e-1  & 1.63 & 8.14 1e-2 & 1.83 \\
1/40&    5.64 1e-1 & 1.06  & 1.23 1e-1 & 1.43 & 4.65 1e-2  & 1.65 & 2.28 1e-2 & 1.83 \\
1/80&    2.49 1e-1 & 1.18  & 4.53 1e-2 & 1.44 & 1.50 1e-2  & 1.64 & 6.47 1e-3 & 1.82 \\
1/160&   1.06 1e-1 & 1.23  & 1.68 1e-2 & 1.43 & 4.90 1e-3  & 1.61 & 1.85 1e-3 & 1.80 \\
1/320&   4.52 1e-2 & 1.24  & 6.30 1e-3 & 1.41 & 1.63 1e-3  & 1.59 & 5.38 1e-4 & 1.79 \\
1/640&   1.92 1e-2 & 1.23  & 2.39 1e-3 & 1.40 & 5.51 1e-4  & 1.57 & 1.58 1e-4 & 1.77 \\
\midrule
\midrule
h & $\alpha=0.9$ & CO & $\alpha=1.2$ & CO & $\alpha=1.5$ & CO & $\alpha=1.8$ & CO \\
\midrule
1/10&    2.16 1e-1 & -     & 1.72 1e-1 & -    & 1.57 1e-1  & -    & 1.53 1e-1  & - \\
1/20&    5.55 1e-2 & 1.96  & 4.17 1e-2 & 2.04 & 3.81 1e-2  & 2.04 & 3.75 1e-2  & 2.03 \\
1/40&    1.42 1e-2 & 1.97  & 1.02 1e-2 & 2.04 & 9.35 1e-3  & 2.03 & 9.28 1e-3  & 2.01 \\
1/80&    3.65 1e-3 & 1.96  & 2.49 1e-3 & 2.03 & 2.31 1e-3  & 2.02 & 2.31 1e-3  & 2.01 \\
1/160&   9.39 1e-4 & 1.96  & 6.11 1e-4 & 2.03 & 5.71 1e-4  & 2.01 & 5.75 1e-4  & 2.00 \\
1/320&   2.42 1e-4 & 1.95  & 1.50 1e-4 & 2.02 & 1.42 1e-4  & 2.01 & 1.44 1e-4  & 2.00 \\
1/640&   6.26 1e-5 & 1.95  & 3.70 1e-5 & 2.01 & 3.53 1e-5  & 2.01 & 3.59 1e-5  & 2.00\\
\bottomrule 
\end{tabular*}\label{table5.15}
\end{table}

\begin{table}[ht]\fontsize{9.5pt}{12pt}\selectfont
\centering \caption{The maximum errors of Improved-Adams methods for (\ref{equa5.11}) when
$t\in[0,1]$.}\vspace{5pt}
\begin{tabular*}{\linewidth}{@{\extracolsep{\fill}}*{9}{c|cc|cc|cc|cc}}
\toprule  
h & $\alpha=0.1$ & CO & $\alpha=0.3$ & CO & $\alpha=0.5$ & CO & $\alpha=0.7$ & CO \\
\midrule
1/10&    2.05 1e+0 & -     & 7.68 1e-1 & -    & 3.74 1e-1  & -    & 2.29 1e-1 & - \\
1/20&    1.13 1e+0 & 0.87  & 2.71 1e-1 & 1.50 & 1.01 1e-1  & 1.89 & 5.32 1e-2 & 2.11 \\
1/40&    5.21 1e-1 & 1.11  & 8.84 1e-2 & 1.61 & 2.59 1e-2  & 1.96 & 1.21 1e-2 & 2.14 \\
1/80&    2.21 1e-1 & 1.24  & 2.81 1e-2 & 1.65 & 6.52 1e-3  & 1.99 & 2.76 1e-3 & 2.13 \\
1/160&   9.04 1e-2 & 1.29  & 8.90 1e-3 & 1.66 & 1.63 1e-3  & 2.00 & 6.37 1e-4 & 2.12 \\
1/320&   3.68 1e-2 & 1.30  & 2.82 1e-3 & 1.66 & 4.05 1e-4  & 2.01 & 1.49 1e-4 & 2.10 \\
1/640&   1.50 1e-2 & 1.30  & 9.00 1e-4 & 1.65 & 1.01 1e-4  & 2.01 & 3.52 1e-5 & 2.08 \\
\midrule
\midrule
h & $\alpha=0.9$ & CO & $\alpha=1.2$ & CO & $\alpha=1.5$ & CO & $\alpha=1.8$ & CO \\
\midrule
1/10&    1.74 1e-1 & -     & 1.50 1e-1 & -    & 1.48 1e-1  & -    & 1.49 1e-1  & - \\
1/20&    3.93 1e-2 & 2.14  & 3.56 1e-2 & 2.08 & 3.61 1e-2  & 2.03 & 3.68 1e-2  & 2.02 \\
1/40&    9.08 1e-3 & 2.11  & 8.69 1e-3 & 2.04 & 8.96 1e-3  & 2.01 & 9.18 1e-3  & 2.00 \\
1/80&    2.15 1e-3 & 2.08  & 2.15 1e-3 & 2.02 & 2.24 1e-3  & 2.00 & 2.29 1e-3  & 2.00 \\
1/160&   5.20 1e-4 & 2.05  & 5.35 1e-4 & 2.01 & 5.58 1e-4  & 2.00 & 5.73 1e-4  & 2.00 \\
1/320&   1.27 1e-4 & 2.03  & 1.34 1e-4 & 2.00 & 1.40 1e-4  & 2.00 & 1.43 1e-4  & 2.00 \\
1/640&   3.15 1e-5 & 2.02  & 3.34 1e-5 & 2.00 & 3.49 1e-5  & 2.00 & 3.58 1e-5  & 2.00\\
\bottomrule 
\end{tabular*}\label{table5.16}
\end{table}

\begin{table}[ht]\fontsize{9.5pt}{12pt}\selectfont
\centering \caption{The maximum errors for (\ref{equa5.11}) when
$\alpha=0.1$, $IN=4$ and $5$.}\vspace{5pt}
\begin{tabular*}{\linewidth}{@{\extracolsep{\fill}}*{5}{c|cc|cc}}  
\toprule  
h & $T=0.1,IN=4$ & CO &$T=0.001,IN=5$ & CO\\
\midrule
T/10&    9.64 1e-9  & -     & 1.45 1e-23 & -     \\
T/20&    6.61 1e-10 & 3.87  & 4.97 1e-25 & 4.87  \\
T/40&    3.91 1e-11 & 4.08  & 1.46 1e-26 & 5.09  \\
T/80&    2.14 1e-12 & 4.19  & 4.01 1e-28 & 5.18  \\
T/160&   1.16 1e-13 & 4.21  & 1.11 1e-29 & 5.18  \\
T/320&   5.83 1e-15 & 4.32  & 2.88 1e-31 & 5.26  \\
T/640&   3.19 1e-16 & 4.19  & 8.19 1e-33 & 5.14  \\
\bottomrule 
\end{tabular*}\label{table5.17}
\end{table}

\begin{figure}
\centering
\includegraphics[width=0.6\textwidth]{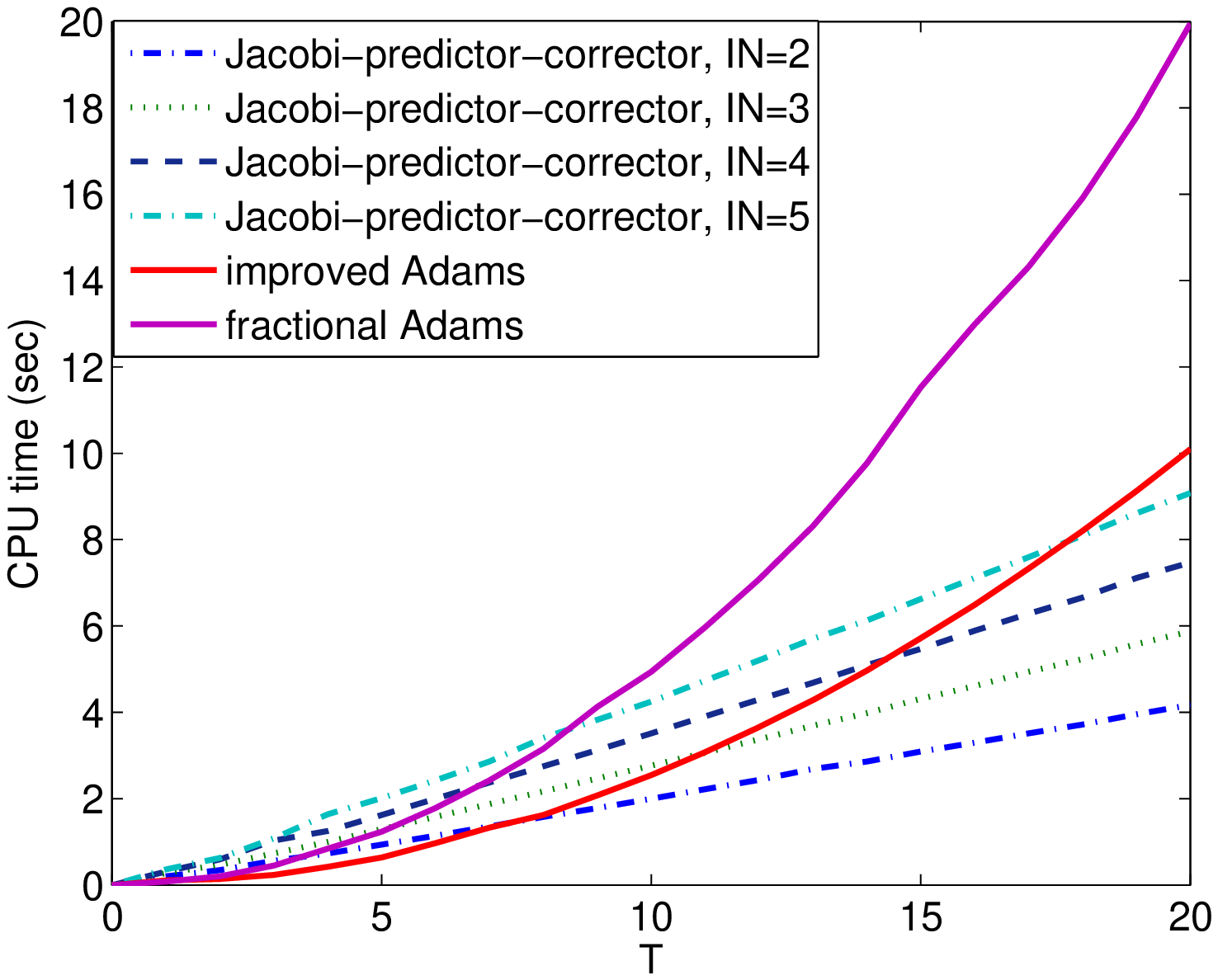}
\caption{The CPU time needed to solve (\ref{equa5.11}) when
$\alpha=0.5,~h=1/40$ for the fractional Adams methods in \cite{Diethelm:02}, the Improved Adams methods
in \cite{Deng:07} and the Jacobi-predictor-corrector approach here when $IN=2,~3,~4,~5$.} \label{fig:511}
\end{figure}

Taking $\alpha=0.5$ and $h=1/40$, a comparison of the CPU time
needed to solve (\ref{equa5.11}) for the fractional Adams methods in
\cite{Diethelm:02}, the Improved Adams methods in \cite{Deng:07} and
the Jacobi-predictor-corrector approach here when $IN=2,~3,~4,~5$,
is reported in Fig. \ref{fig:511}. Fig. \ref{fig:511} illustrates
that the computational cost of the Jacobi-predictor-corrector
approach is $O(N)$. Also notice that, as expected, both the
fractional Adams methods and the Improved Adams methods exhibit
$O(N^2)$ computational complexity.

Although, from Fig. \ref{fig:511}, the Jacobi-predictor-corrector approach takes more time to reach the
terminate time T, when T is small,
for example,
when $T=1$, the CPU time of the fractional Adams methods, the Improved Adams methods and
the Jacobi-predictor-corrector approach when $IN=2,~3,~4,~5$, respectively are
$0.078,~0.109,~0.203,~0.281,~0.313,$ and $0.359$ (sec). While, by Table \ref{table5.11} to Table \ref{table5.16},
we can see that, when $N=160$, the maximum errors of the Jacobi-predictor-corrector methods,
$9.99$~1e-3,~$7.19$~1e-4,~$1.43$~1e-5,~$1.08$~1e-6, are much smaller than those of the other two methods, $4.65$~1e-2,~$2.95$~1e-2.

Table \ref{table5.18} shows the CPU time (sec) and the steps $N$
needed to solve (\ref{equa5.11}) when $\alpha=0.5$ with the maximum
error $1.0\times10^{-3}$, for the fractional Adams methods in
\cite{Diethelm:02}, the Improved Adams methods in \cite{Deng:07} and
the Jacobi-predictor-corrector approach here when $IN=2,~3,~4,~5$.
The consumed CPU time presented in Table \ref{table5.18} shows that
the fractional Adams methods generates the numerical solution with
the same accuracy as the other two methods, but uses much less CPU
time. This advantage is more obvious as the terminate time goes
long. It further demonstrates the quickness and efficiency of the
Jacobi-predictor-corrector method.
\begin{table}[ht]\fontsize{9.5pt}{12pt}\selectfont
\centering \caption{The CPU time (sec) and the steps $N$ needed to
solve (\ref{equa5.11}) when $\alpha=0.5$ with the maximum error
$1.0\times10^{-3}$, for the fractional Adams methods in
\cite{Diethelm:02}, the Improved Adams methods in \cite{Deng:07} and
the Jacobi-predictor-corrector approach here when
$IN=2,~3,~4,~5$.}\vspace{5pt}
\begin{tabular*}{\linewidth}{@{\extracolsep{\fill}}*{10}{c|c|cc|cc|cc|cc}} 
\toprule  
\multicolumn{2}{c|}{}&\multicolumn{8}{c}{terminal time}\\
\cline{3-10}
\multicolumn{2}{c|}{methods}&\multicolumn{2}{c|}{$T=0.5$}&\multicolumn{2}{c|}{$T=1.0$}&\multicolumn{2}{c}{$T=1.5$}&\multicolumn{2}{c}{$T=2.0$}\\
\cline{3-10}
\multicolumn{2}{c|}{}&$N$&CPU time (sec)&$N$&CPU time (sec)&$N$&CPU time (sec)&$N$&CPU time (sec)\\
\midrule
\multicolumn{2}{c|}{fractional Adams}&17  &6.25 1e -2 &432   &5.72 1e+0   &3240 &3.20 1e+2  &14200 &6.26 1e+3\\
\multicolumn{2}{c|}{Improved Adams}  &14  &3.13 1e -2 &204   &7.97 1e -1  &945  &1.41 1e+1  &2831  &1.29 1e+2\\
\cline{1-2}
Jacobi-   &$IN=2$                    &11  &9.34 1e -2 &119   &6.09 1e -1  &492  &2.52 1e+0  &1456  &7.75 1e+0  \\
predictor-&$IN=3$                    &7   &1.56 1e -2 &34    &1.41 1e -1  &89   &5.16 1e -1 &117   &1.13 1e+0\\
corrector &$IN=4$                    &5   &1.56 1e -2 &18    &7.81 1e -2  &34   &1.72 1e -1 &51    &2.97 1e -1\\
methods   &$IN=5$                    &5   &3.13 1e -2 &13    &4.69 1e -2  &23   &1.09 1e -1 &33    &1.88 1e -1\\
\bottomrule 
\end{tabular*}\label{table5.18}
\end{table}

\subsection{Robustness of the algorithm}\label{sec:5.2}
Here we study the following equation as an example to show the
robustness of the algorithm,
\begin{equation}\label{equa5.21}
D^\alpha_\ast x(t)=-x(t),~~~~~~x(0)=1,~~~x'(0)=0 ~~(\textrm{if} ~1<\alpha\leq 2).
\end{equation}
It is well known that the exact solution of (\ref{equa5.21}) is
\begin{equation}\label{equa5.22}
x(t)=E_{\alpha}(- t^{\alpha}),
\end{equation}
where
\begin{equation}\label{equa5.23}
E_{\alpha}(z)=\sum_{k=0}^{\infty}\frac{z^k}{\Gamma(\alpha k +1)}, ~~~Re(\alpha)>0,z\in \mathbb{C},
\end{equation}
is the Mittag-Leffler function of order $\alpha$. It is obvious that
neither $x(t)$ nor $D^\alpha_\ast x(t)$ has a bounded first (second)
derivative at $t=0$ when $0<\alpha \leq 1$ ($1<\alpha \leq 2$), so
we deal with (\ref{equa5.21}) as we discussed in Section
\ref{sec:4.1}. Here we take $T_0=0.1$, $JN=26$, $\tilde{JN}=2JN$,
where $T_0,~JN$ and $\tilde{JN}$ are defined as in Section
\ref{sec:4.1}, and the exact solutions are calculated using the
function $``mlf.m"$ \cite{Podlubny}. The convergent order is also
simply verified in Table \ref{table5.21} and Table \ref{table5.22}.

\begin{table}[ht]\fontsize{9.5pt}{12pt}\selectfont
\centering \caption{The maximum errors for (\ref{equa5.21}) when
$t\in[0,1.1]$ and $IN=2$.}\vspace{5pt}
\begin{tabular*}{\linewidth}{@{\extracolsep{\fill}}*{9}{c|cc|cc|cc|cc}}
\toprule  
h & $\alpha=0.2$ & CO & $\alpha=0.5$ & CO & $\alpha=1.2$ & CO & $\alpha=1.8$ & CO \\
\midrule
1/10&    4.84 1e-3 & -     & 2.30 1e-3 & -    & 1.20 1e-4  & -    & 3.72 1e-4 & - \\
1/20&    1.49 1e-3 & 1.70  & 5.10 1e-4 & 2.17 & 3.47 1e-5  & 1.79 & 1.06 1e-4 & 1.81 \\
1/40&    4.04 1e-4 & 1.88  & 1.02 1e-4 & 2.33 & 7.83 1e-6  & 2.15 & 2.62 1e-5 & 2.02 \\
1/80&    9.96 1e-5 & 2.02  & 1.89 1e-5 & 2.43 & 1.87 1e-6  & 2.07 & 6.35 1e-6 & 2.05 \\
1/160&   2.44 1e-5 & 2.03  & 3.95 1e-6 & 2.26 & 5.41 1e-7  & 1.79 & 1.62 1e-6 & 1.97 \\
\bottomrule 
\end{tabular*}\label{table5.21}
\end{table}

\begin{table}[ht]\fontsize{9.5pt}{12pt}\selectfont
\centering \caption{The maximum errors for (\ref{equa5.21}) when
$t\in[0,1.1]$ and $IN=3$.}\vspace{5pt}
\begin{tabular*}{\linewidth}{@{\extracolsep{\fill}}*{9}{c|cc|cc|cc|cc}}
\toprule  
h & $\alpha=0.2$ & CO & $\alpha=0.5$ & CO & $\alpha=1.2$ & CO & $\alpha=1.8$ & CO \\
\midrule
1/10&    2.77 1e-3 & -     & 7.30 1e-4 & -    & 1.11 1e-5  & -    & 1.64 1e-5 & - \\
1/20&    6.04 1e-4 & 2.20  & 1.14 1e-4 & 2.68 & 3.23 1e-6  & 1.78 & 3.00 1e-6 & 2.45 \\
1/40&    1.06 1e-4 & 2.51  & 1.43 1e-5 & 2.99 & 5.48 1e-7  & 2.56 & 4.64 1e-7 & 2.69 \\
1/80&    1.29 1e-5 & 3.04  & 1.40 1e-6 & 3.36 & 8.88 1e-8  & 2.63 & 5.91 1e-8 & 2.97 \\
1/160&   1.36 1e-6 & 3.25  & 3.78 1e-8 & 5.21 & 1.09 1e-8  & 3.03 & 7.84 1e-9 & 2.92 \\
\bottomrule 
\end{tabular*}\label{table5.22}
\end{table}

\begin{figure}[!htb]
    \subfigure[]{
    \begin{minipage}[t]{0.48\linewidth}
        \centering
        \includegraphics[scale=0.3]{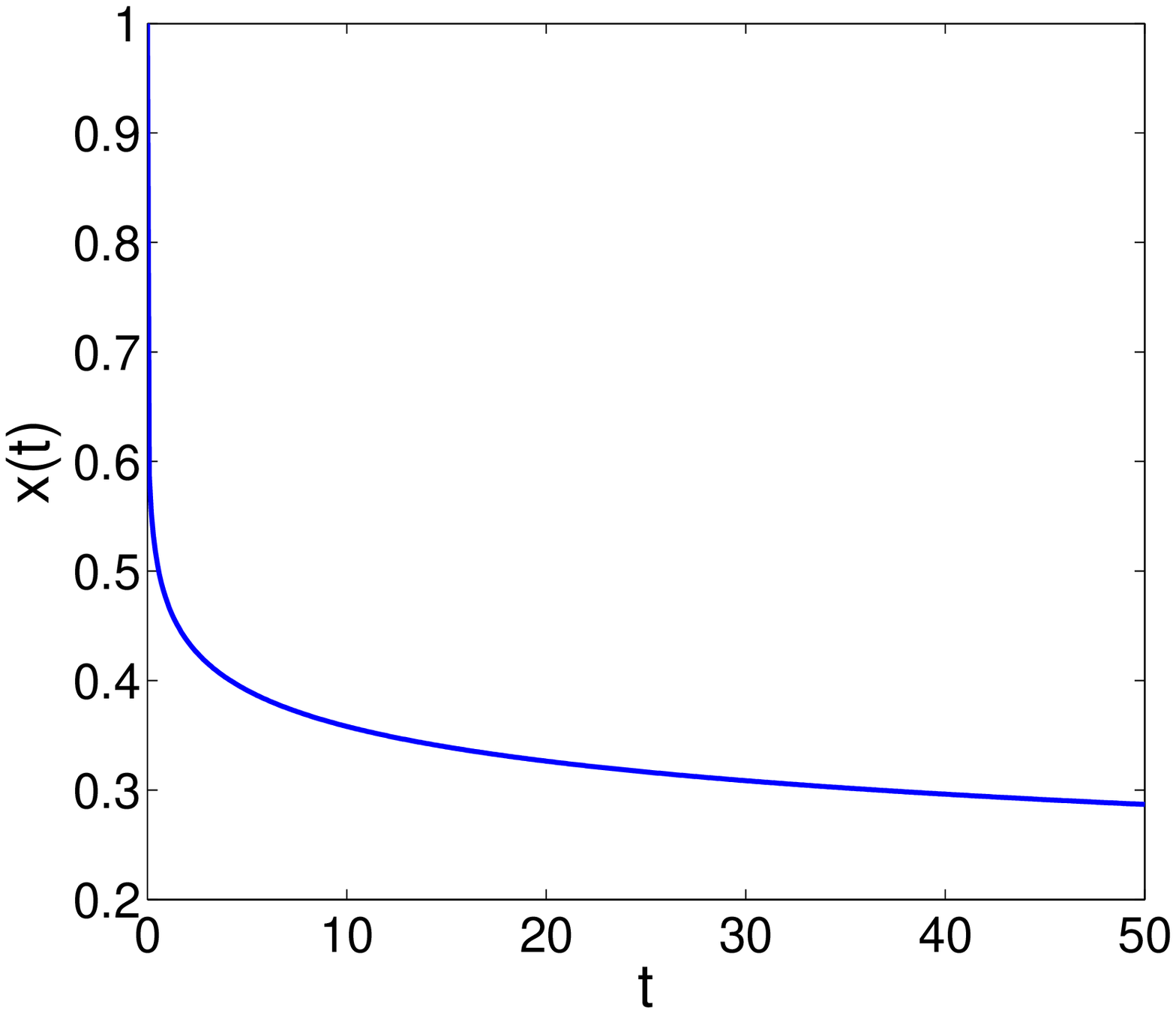}
      \end{minipage}}
      \subfigure[]{
      \begin{minipage}[t]{0.48\linewidth}
        \centering
        \includegraphics[scale=0.3]{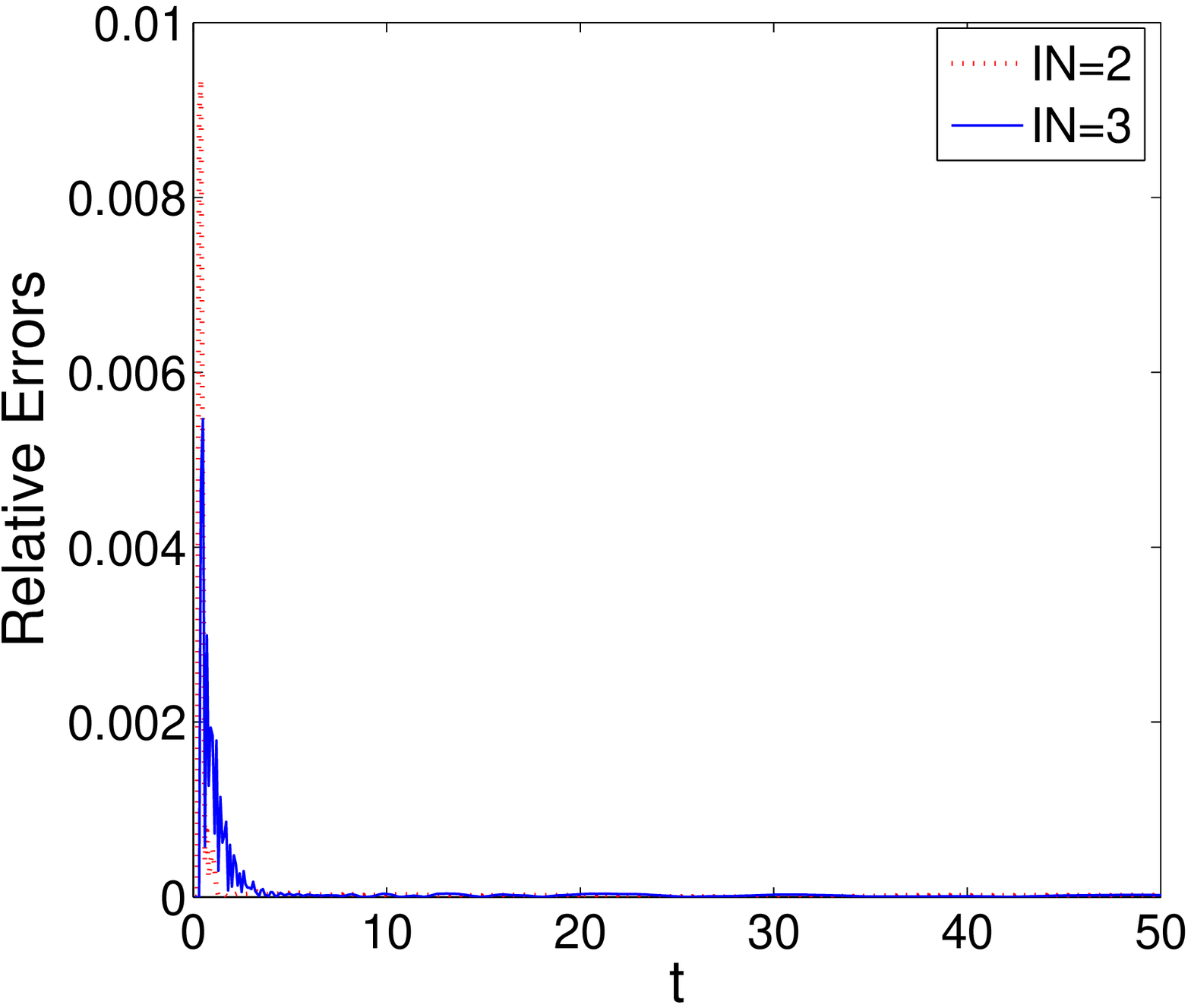}
      \end{minipage}}
      \caption{(a) Exact solution of (\ref{equa5.21}) for $\alpha=0.2$;
      (b) Relative errors with $\alpha = 0.2$, $h=1/10$, $IN=2$ or $3$.}\label{fig:521}
    \subfigure[]{
    \begin{minipage}[t]{0.48\linewidth}
        \centering
        \includegraphics[scale=0.3]{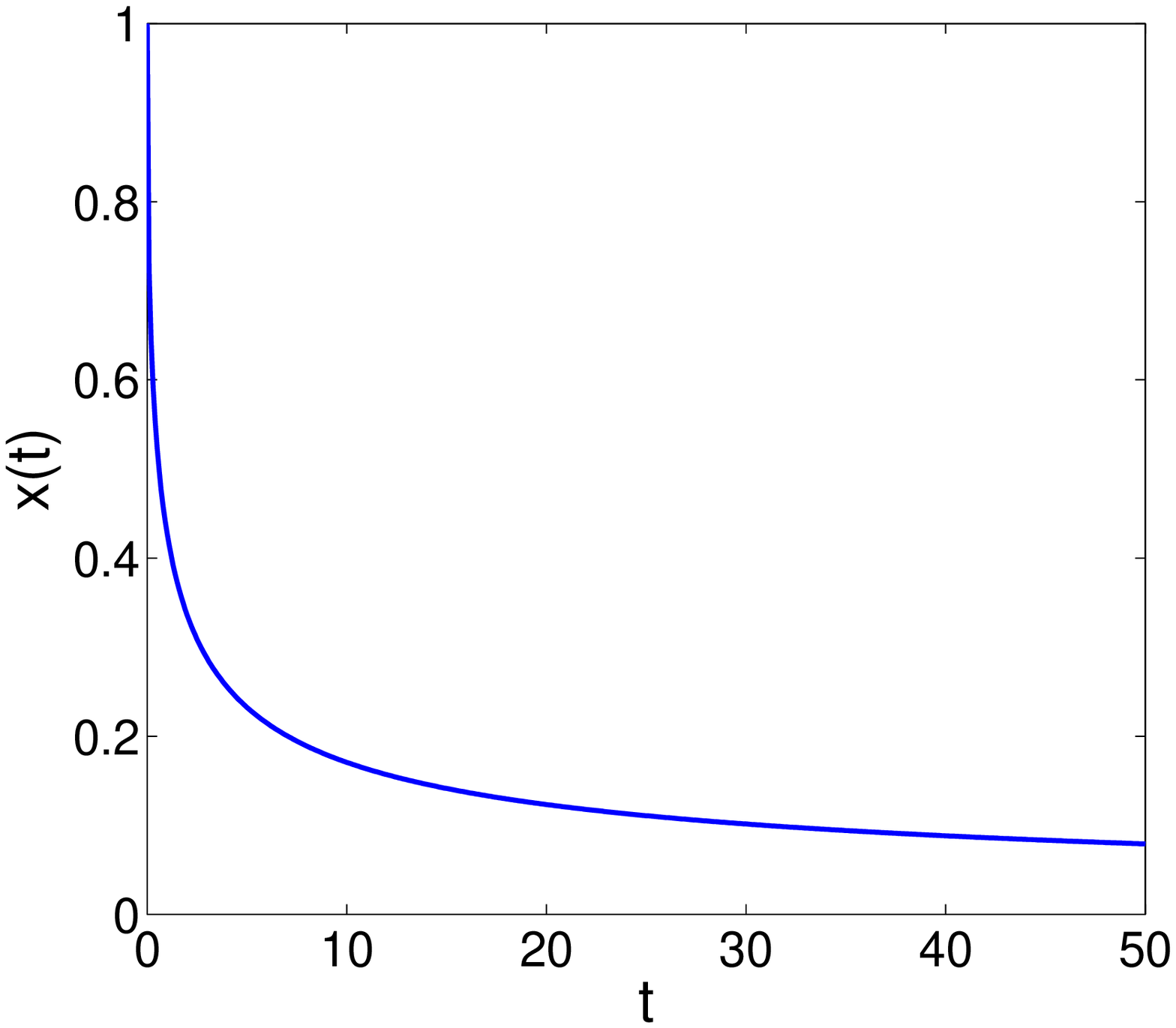}
      \end{minipage}}
      \subfigure[]{
      \begin{minipage}[t]{0.48\linewidth}
        \centering
        \includegraphics[scale=0.3]{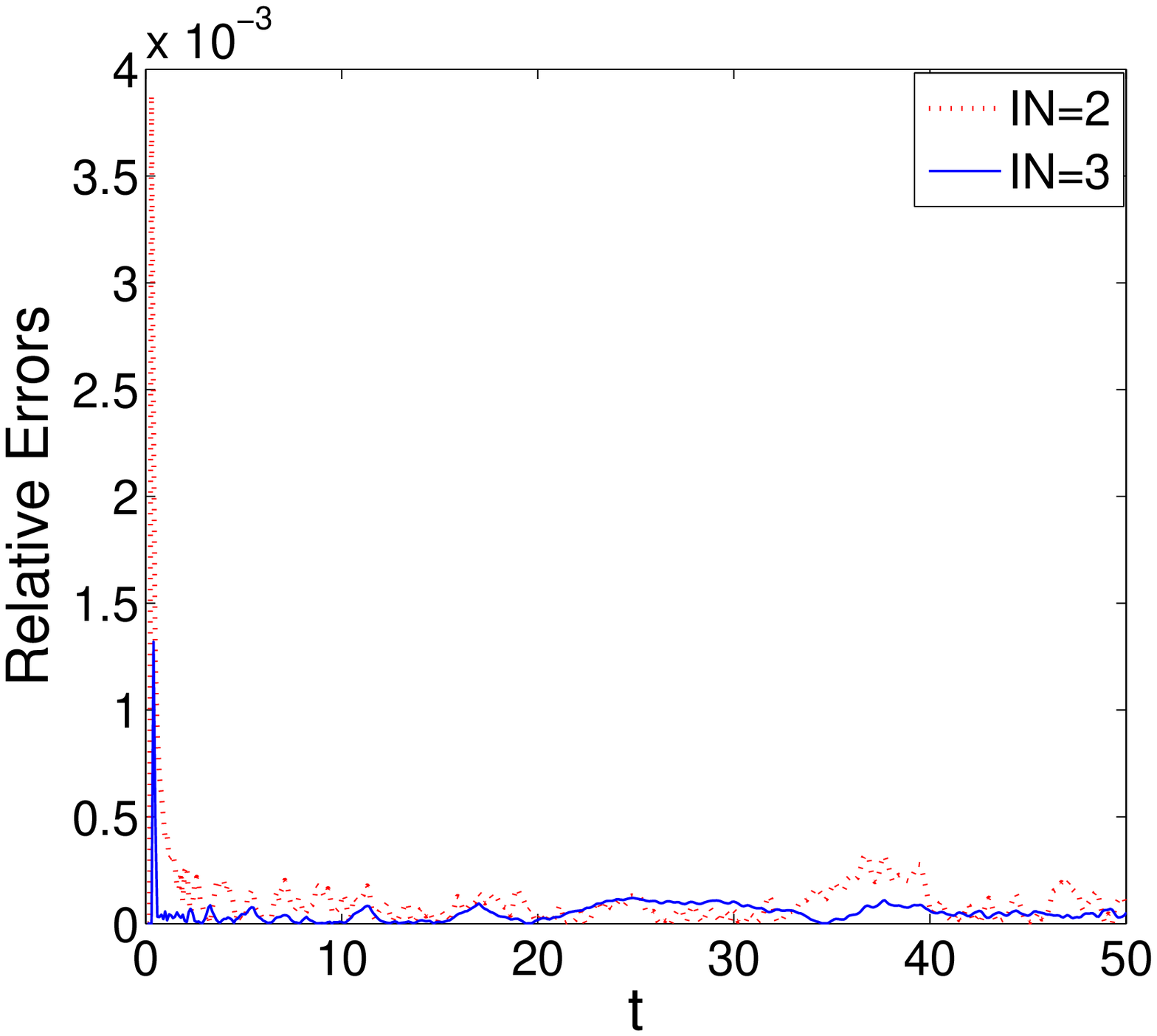}
      \end{minipage}}
      \caption{(a) Exact solution of (\ref{equa5.21}) for $\alpha=0.5$;
      (b) Relative errors with $\alpha = 0.5$, $h=1/10$, $IN=2$ or $3$.}\label{fig:522}
\end{figure}

Further we compute (\ref{equa5.21}) with a big time interval,
$T=50$; the parameters $IN,~JN,~\tilde{JN}$ are taken the same as
the above ones and $h$ is taken as $1/10$. The exact solutions and
relative errors are shown in Fig. \ref{fig:521} with $\alpha=0.2$ and Fig. \ref{fig:522}
with $\alpha=0.5$. It can be seen that the relative errors in the
interval are less than $O(10^{-4})$ when time is long, which
suggests that our method is suitable for the long-time computation.

\section{Conclusion}\label{sec:6}
We provide the Jacobi-predictor-corrector approach for the fractional ordinary differential
equations; the basic idea is to take the Riemann-Liouville integral kernel as the Jacobi-weight
function, and to realize the algorithm by doing the Jacobi-Gauss-Lobatto quadrature and polynomial
interpolation. The convergent order is exactly equal to the number of interpolating nodes $IN$. The
computational complexity is $O(N)$ for $\alpha \in (0,\infty)$, where $N$ is the total
computational steps. This is the striking feature/advantage of the algorithm, since the
computational complexity of numerically solving the fractional ordinary differential equation
usually is $O(N^2)$, caused by its nonlocal property; when $\alpha \in (0,2)$, it is possible to
reduce the computational cost to $O(N\log N)$ by combining the short memory principle.

\bigskip
\noindent {\bf Acknowledgments.} This research was supported by the Program for New Century Excellent
Talents in University under Grant No. NCET-09-0438, the National
Natural Science Foundation of China under Grant No. 10801067, and
the Fundamental Research Funds for the Central Universities under
Grant No. lzujbky-2010-63 and No. lzujbky-2012-k26. WHD thanks
Chi-Wang Shu for the discussions and valuable comments.

\clearpage


\section*{Appendix}

We display the Jacobi-Gauss-Lobatoo nodes and weighs in the
reference interval $[-1,1]$ used in the numerical experiments in the
following tables, where the weight function is $(1-s)^\alpha$,
$\alpha=0.1$, $0.3$, $0.5$, $0.7$, $0.9$, $1.2$, $1.5$, $1.8$,
respectively, and the number of the quadrature nodes $JN+1=27$.

\begin{table}[ht]\fontsize{9.5pt}{12pt}\selectfont
\begin{center} \caption{}\vspace{5pt}
\begin{tabular*}{\linewidth}{@{\extracolsep{\fill}}*{9}{cc|cc}}
\toprule  
\multicolumn{2}{c|}{~$\alpha=0.1$}
&\multicolumn{2}{c}{~$\alpha=0.3$}\\
\midrule
nodes & weights &nodes & weights \\
\midrule
-1.0000000000000000 &0.0015793891284060 &-1.0000000000000000 &0.0018004451789191\\
-0.9892016529048960 &0.0097486968513111 &-0.9892836220495154 &0.0111011132499144\\
-0.9639539701336150 &0.0176099210456689 &-0.9642264317182587 &0.0200016826109024\\
-0.9247048615099546 &0.0255586083473313 &-0.9252701785814652 &0.0289128339852526\\
-0.8720275885827603 &0.0336423948344829 &-0.8729795316392537 &0.0378465851188033\\
-0.8066879161228238 &0.0419097843832046 &-0.8081088436911802 &0.0468128223976810\\
-0.7296351534780572 &0.0504145465068031 &-0.7315934254809393 &0.0558230223195378\\
-0.6419886593191621 &0.0592179803160306 &-0.6445363551226361 &0.0648912240850699\\
-0.5450216478217176 &0.0683915067416097 &-0.5481926424422340 &0.0740349165644332\\
-0.4401427133803333 &0.0780200133502234 &-0.4439511569824626 &0.0832761166251512\\
-0.3288753760288977 &0.0882063337718825 &-0.3333146137264882 &0.0926427678850384\\
-0.2128359531749780 &0.0990774042058216 &-0.2178779134712215 &0.1021706112368920\\
-0.0937100816930866 &0.1107929333210531 &-0.0993051526697154 &0.1119057485004456\\
0.0267717676524930 &0.1235579229632447   &0.0206943649229667 &0.1219082421622975\\
0.1468594274088615 &0.1376412507907670   &0.1403907686322420 &0.1322573016048666\\
0.2648084569648291 &0.1534041057504351   &0.2580585579193416 &0.1430589703422321\\
0.3789054830602389 &0.1713450510950583   &0.3720014765729532 &0.1544578938601842\\
0.4874930892080245 &0.1921744216247991   &0.4805769654678000 &0.1666560217001116\\
0.5889938926263092 &0.2169432942460124   &0.5822198414675918 &0.1799436799255219\\
0.6819334593243178 &0.2472807617819305   &0.6754648613490701 &0.1947540530811646\\
0.7649617255263299 &0.2858641340356679   &0.7589678460750942 &0.2117653286898140\\
0.8368726174182031 &0.3374445171968009   &0.8315250626755174 &0.2321093791361104\\
0.8966215953237894 &0.4113924209918509   &0.8920905904778630 &0.2578499674265770\\
0.9433409167504676 &0.5293048192507985   &0.9397914478954834 &0.2932722023123962\\
0.9763527097958942 &0.7555606909447271   &0.9739404234643343 &0.3493613894433857\\
0.9951835446365114 &1.4176717079984273   &0.9940482649436114 &0.4687789976411144\\
1.0000000000000000 &5.0539800138885518   &1.0000000000000000 &0.4664213940659055\\
\bottomrule 
\end{tabular*}\label{tablesp.1}
\end{center}
\end{table}
\clearpage
\begin{table}[ht]\fontsize{9.5pt}{12pt}\selectfont
\begin{center} \caption{}\vspace{5pt}
\begin{tabular*}{\linewidth}{@{\extracolsep{\fill}}*{9}{cc|cc}}
\toprule  
\multicolumn{2}{c|}{~$\alpha=0.5$}
&\multicolumn{2}{c}{~$\alpha=0.7$}\\
\midrule
nodes & weights &nodes & weights \\
\midrule
-1.0000000000000000 &0.0020525595970582 &-1.0000000000000000 &0.0023401106204443\\
-0.9893643555933919 &0.0126420869565904 &-0.9894438812775599 &0.0143980191213123\\
-0.9644947993303408 &0.0227208651492852 &-0.9647591646875335 &0.0258125839599000\\
-0.9258270443312571 &0.0327129904967510 &-0.9263756474471649 &0.0370189772483234\\
-0.8739173433922172 &0.0425872817104095 &-0.8748413378464713 &0.0479340538950119\\
-0.8095088697871018 &0.0523089398061103 &-0.8108884559986644 &0.0584716242404765\\
-0.7335232176624819 &0.0618432836433182 &-0.7354251541586265 &0.0685470018011806\\
-0.6470475078824108 &0.0711562189729229 &-0.6495229109880948 &0.0780777884066716\\
-0.5513189009359872 &0.0802144213662633 &-0.5544013837390309 &0.0869842827280500\\
-0.4477069175659654 &0.0889854705911301 &-0.4514111111711064 &0.0951897677589023\\
-0.3376938533233919 &0.0974379713316180 &-0.3420143458160455 &0.1026206972731357\\
-0.2228535755897699 &0.1055416672643979 &-0.2277642960175138 &0.1092067627458629\\
-0.1048290089639729 &0.1132675500578416& -0.1102830750128124 &0.1148808056093876\\
0.0146913679757165 &0.1205879635232249 &0.0087613289573869 &0.1195785218117537\\
0.1339976779295242 &0.1274767027659907 &0.1276787341676305 &0.1232378787845395\\
0.2513831064264158 &0.1339091080692547 &0.2447807621576336 &0.1257981209039414\\
0.3651683196413840 &0.1398621532108058 &0.3584048089141069 &0.1271981640581311\\
0.4737254891175540 &0.1453145279145870 &0.4669376500627146 &0.1273740449364915\\
0.5755015796722314 &0.1502467141498365 &0.5688383448113328 &0.1262548372106364\\
0.6690405673158345 &0.1546410560090228 &0.6626601132252578 &0.1237559426680637\\
0.7530042693246988 &0.1584818229166975 &0.7470708756553147 &0.1197675894104367\\
0.8261914884659855 &0.1617552659442731 &0.8208721610920136 &0.1141338609152897\\
0.8875551974923557 &0.1644496670298319 &0.8830161106263974 &0.1066110258518362\\
0.9362175180611009 &0.1665553809272196 &0.9326203125145781 &0.0967739548306226\\
0.9714822797862176 &0.1680648697345029 &0.9689801330973330 &0.0837631042002331\\
0.9928449797512120 &0.1689727298783471 &0.9915771271381129 &0.0653397616047768\\
1.0000000000000000 &0.0846378557289007 &1.0000000000000000 &0.0196518498509760\\
\bottomrule 
\end{tabular*}\label{tablesp.2}
\end{center}
\end{table}
\clearpage

\begin{table}[ht]\fontsize{9.5pt}{12pt}\selectfont
\begin{center} \caption{}\vspace{5pt}
\begin{tabular*}{\linewidth}{@{\extracolsep{\fill}}*{9}{cc|cc}}
\toprule  
\multicolumn{2}{c|}{~$\alpha=0.9$}
&\multicolumn{2}{c}{~$\alpha=1.2$}\\
\midrule
nodes & weights &nodes & weights \\
\midrule
-1.0000000000000000 &0.0026680954862362 &-1.0000000000000000 &0.0032485813206930\\
-0.9895222260184334 &0.0163990211169060 &-0.9896375857934652 &0.0199364618399667\\
-0.9650196167842330 &0.0293282050572500 &-0.9654031458102729 &0.0355265819180070\\
-0.9269161710244489 &0.0418988082171859 &-0.9277121945697792 &0.0504662938531113\\
-0.8757518197224115 &0.0539656088052377 &-0.8770928511031259 &0.0644948884276532\\
-0.8122480503073657 &0.0653836420874911 &-0.8142509047810206 &0.0773638096423790\\
-0.7372998407745818 &0.0760144717348032 &-0.7400620590862488 &0.0888474771253680\\
-0.6519633346267486 &0.0857281246221285 &-0.6555600135812010 &0.0987484771441300\\
-0.5574410233343691 &0.0944045439109477 &-0.5619221255221779 &0.1069017533569795\\
-0.4550648215271042 &0.1019348490180515 &-0.4604530253538590 &0.1131780974268036\\
-0.3462773061846129 &0.1082224273322169 &-0.3525664467804354 &0.1174869327523327\\
-0.2326113924086162 &0.1131838455110155 &-0.2397655327004557 &0.1197783563983938\\
-0.1156687346912121 &0.1167495593531492 &-0.1236218938939909 &0.1200444135582023\\
0.0029028410706502 &0.1188643960701240 &-0.0057537131214067 &0.1183195966655319\\
0.1214325564431927 &0.1194877759976240 &0.1121967999388062 &0.1146805830124639\\
0.2382502228570740 &0.1185936294194118 &0.2285862886357495 &0.1092452502518708\\
0.3517097756999029 &0.1161699441834089 &0.3417931453037846 &0.1021710400770435\\
0.4602124683760809 &0.1122178437862630 &0.4502401041519408 &0.0936527803497127\\
0.5622293993378423 &0.1067500287515353 &0.5524162161956613 &0.0839201325248429\\
0.6563230543010921 &0.0997882844987957 &0.6468978998418504 &0.0732349201079911\\
0.7411675592086939 &0.0913594924131627 &0.7323687733098204 &0.0618887495635877\\
0.8155673560609422 &0.0814889890869870 &0.8076379911244389 &0.0502016389583398\\
0.8784740302618931 &0.0701886634133376 &0.8716568246106763 &0.0385230330239028\\
0.9290010193400755 &0.0574330648431000 &0.9235332376746080 &0.0272382282726720\\
0.9664358148554766 &0.0431024991510826 &0.9625441648076186 &0.0167880196253709\\
0.9902477964312976 &0.0268013946767667 &0.9881446058894564 &0.0077265395634170\\
1.0000000000000000 &0.0052794393153508 &1.0000000000000000 &0.0008846215676251\\
\bottomrule 
\end{tabular*}\label{tablesp.3}
\end{center}
\end{table}
\clearpage
\begin{table}[ht]\fontsize{9.5pt}{12pt}\selectfont
\begin{center} \caption{}\vspace{5pt}
\begin{tabular*}{\linewidth}{@{\extracolsep{\fill}}*{9}{cc|cc}}
\toprule  
\multicolumn{2}{c|}{~$\alpha=1.5$}
&\multicolumn{2}{c}{~$\alpha=1.8$}\\
\midrule
nodes & weights &nodes & weights \\
\midrule
-1.0000000000000000 &0.0039558421325122 &-1.0000000000000000 &0.0048176566867522\\
-0.9897504320058830 &0.0242407262612330 &-0.9898608459541634 &0.0294787408291468\\
-0.9657783446854108 &0.0430450442150069 &-0.9661454822450298 &0.0521665283985357\\
-0.9284910158577785 &0.0608075612568833 &-0.9292531881022244 &0.0732935305513724\\
-0.8784051049048023 &0.0771197426976410 &-0.8796895020922768 &0.0922636940579904\\
-0.8162111679198959 &0.0916086237361090 &-0.8181301941668921 &0.1085561087536052\\
-0.7427661984333459 &0.1039548961583067 &-0.7454140912397981 &0.1217534509743617\\
-0.6590821003162865 &0.1139030192180424 &-0.6625319256813679 &0.1315582543475505\\
-0.5663118098346140 &0.1212688126551674 &-0.5706128999526579 &0.1378033510671248\\
-0.4657334305951538 &0.1259447610428847 &-0.4709093210906420 &0.1404566089311143\\
-0.3587326320143288 &0.1279029357867882 &-0.3647795460854125 &0.1396198574132639\\
-0.2467835618430933 &0.1271954603677246 &-0.2536694784396956 &0.1355220825893431\\
-0.1314285381791132 &0.1239525097131300 &-0.1390928704191656 &0.1285072066852624\\
-0.0142568016179130 &0.1183779081926168 &-0.0226107001846948 &0.1190169991802236\\
0.1031173793924301 &0.1107424642216554 &0.0941900949493554 &0.1075698748132034\\
0.2190769482282871 &0.1013752491912932 &0.2097182393627444 &0.0947365082759497\\
0.3320243383882507 &0.0906530919063302 &0.3223997996747073 &0.0811133270750175\\
0.4404034835594873 &0.0789886148348753 &0.4306996300173916 &0.0672950260065822\\
0.5427212564354959 &0.0668171834999999 &0.5331422907681873 &0.0538472733927563\\
0.6375680414418941 &0.0545831738546201 &0.6283321577931673 &0.0412807461037388\\
0.7236371592440043 &0.0427259833787319 &0.7149724532569940 &0.0300275324422679\\
0.7997428790488688 &0.0316662190978989 &0.7918829524158366 &0.0204207704286011\\
0.8648367821754457 &0.0217924888159768 &0.8580161681069255 &0.0126781227682169\\
0.9180222971008329 &0.0134491973972607 &0.9124719559264060 &0.0068892709351844\\
0.9585674487321285 &0.0069256857723223 &0.9545111721119047 &0.0030068635217618\\
0.9859178863652559 &0.0024466760492113 &0.9835752524825214 &0.0008386286076415\\
1.0000000000000000 &0.0001742117099051 &1.0000000000000000 &0.0000387924881512\\
\bottomrule 
\end{tabular*}\label{tablesp.4}
\end{center}
\end{table}

\clearpage


\end{document}